\documentclass{amsart}
\usepackage{amssymb, graphics}
\usepackage[all]{xy}
\newcommand\datver[1]{\def\datverp%
 {\par\boxed{\boxed{\text{#1; Run: \today}}}}}

\numberwithin{equation}{section}
\newtheorem{lemma}{Lemma}[section]
\newtheorem{proposition}[lemma]{Proposition}
\newtheorem{corollary}[lemma]{Corollary}
\newtheorem{transversality}[lemma]{Transversality theorem}
\newtheorem{extension}[lemma]{Extension theorem}
\newtheorem{theorem}{Theorem}
\newtheorem{non-theorem}[lemma]{Non-Theorem}

\newtheorem{definition}[lemma]{Definition}
\newtheorem{remark}[lemma]{Remark}

\newcommand\cf{cf\@. }
\newcommand\ie{i\@.e\@. }
\newcommand\eg{e\@.g\@. }

\newcommand\pa{ \partial}

\newcommand\bbC{\mathbb C}
\newcommand\bbD{\mathbb D}

\newcommand\bbN{\mathbb N}
\newcommand\bbP{\mathbb P}

\newcommand\bbR{\mathbb R}
\newcommand\bbS{\mathbb S}

\newcommand\bbZ{\mathbb Z}

\newcommand\cQ{\mathcal{Q}}

\newcommand\cU{\mathcal{U}}
\newcommand\tcU{\widetilde{\mathcal{U}}}
\newcommand\cV{\mathcal{V}}
\newcommand\cZ{\mathcal{Z}}

\newcommand\tV{\widetilde{V}}
\newcommand\tS{\widetilde{S}}
\newcommand\cW{\mathcal{W}}

\newcommand\End{\operatorname{End}}

\newcommand\nb{ N  }   
\newcommand\tnb[1]{\widetilde{\nb{#1}}}
\newcommand\unb{\bbS N }  
\newcommand\spa{\operatorname{span}}
\newcommand\tM{\widetilde{M}}
\newcommand\tI{\widetilde{I}}
\newcommand\bI{\overline{I}}
\newcommand\bS{\overline{S}}
\newcommand\hI{\hat{I}}
\newcommand\hp{\hat{p}}
\newcommand\tp{\widetilde{p}}
\newcommand\tD{\widetilde{D}}
\newcommand\hD{\hat{D}}
\newcommand\bD{\overline{D}}

\newcommand\tZ{\widetilde{Z}}

\newcommand\tnu{\widetilde{\nu}}
\newcommand\ho{\operatorname{H}}
\newcommand\CI{ \mathcal{C}^{\infty}}

\newcommand\hcf{ h_{\phi}}  
\newcommand\tSigma{\widetilde{\Sigma}}

\newcommand\cp[1]{\bbC\bbP_{#1}}  
\newcommand\tcp[1]{\widetilde{\bbC\bbP}_{#1}}  
\newcommand\rp[1]{\bbR\bbP^{#1}}  

\newcommand\SO{\operatorname{SO}}
\newcommand\Span{\operatorname{span}}

\datver{0.2A; Revised: May 21, 2008}
\begin{document}
\title[Families of holomorphic disks]
{On the uniqueness of certain families of holomorphic disks}

\author{Fr\'ed\'eric Rochon}
\address{Department of mathematics, University of Toronto}
\email{rochon@math.utoronto.ca}
\dedicatory{\datverp}
\thanks{The author aknowledge the support of the Fonds qu\'eb\'ecois de la recherche sur la nature et les
technologies while part of this work was conducted.}

\begin{abstract}   
A Zoll metric is a Riemannian metric whose geodesics are all circles of equal length.  Via the twistor correspondence of LeBrun and Mason, a Zoll metric on the sphere
$\bbS^{2}$ corresponds to a family of holomorphic disks in $\cp{2}$ with boundary in a 
totally real submanifold $P\subset\cp{2}$.  In this paper, we show that for a fixed $P\subset \cp{2}$, such
a family is unique if it exists, implying that the twistor correspondence of LeBrun and Mason is injective.  One of the key ingredients in the proof is the blow-up and blow-down constructions in the sense of Melrose. 
\end{abstract}
\maketitle


\section*{Introduction}

Since their introduction by Roger Penrose \cite{Penrose} (see also \cite{Atiyah-Hitchin-Singer}), twistor spaces have proven to be very successful in the study of Riemannian manifolds with special properties.  They
allow one to encode subtle geometric information into holomorphic objects.  In this way, powerful tools
of complex and algebraic geometry can come into play and reveal in a rather unexpected way new properties
of the original geometric situation.  A recent example of such a phenomenon was discovered by 
LeBrun and Mason \cite{Lebrun-Mason} in their study of Zoll metrics on a compact surface.

In general, a Zoll metric on a manifold is a metric whose geodesics are simple closed curves of equal
length.  The terminology is in honor of Otto Zoll \cite{Zoll}, who discovered that the 2-dimensional 
sphere $\bbS^{2}$ has many such metrics beside the standard one.  Later on, it was speculated by Funk
\cite{Funk}
and proven by Guillemin \cite{Guillemin_Zoll} that the tangent space of the moduli space of Zoll metrics on $\bbS^{2}$ (modulo isometries and rescaling) at the standard round metric is isomorphic to the space of 
odd functions $f:\bbS^{2}\to \bbR$.  In particular, this indicates that the moduli space of Zoll metrics
on $\bbS^{2}$ is infinite dimensional.  The only other compact surface admitting a Zoll metric, 
$\rp{2}$, has a very different behavior in that respect since, as was conjectured by Blaschke and 
proved by Leon Green \cite{Green}, it admits only one Zoll metric modulo rescaling and isometries. 

Via the use of twistor theory, LeBrun and Mason were able in \cite{Lebrun-Mason} to recover all these results in a very elegant
and novel way.  In the case of the sphere $\bbS^{2}$, their twistor correspondence associates to a Zoll
metric on $\bbS^{2}$ a family of holomorphic disks in $\cp{2}$ with boundary in a totally real submanifold
$P\subset \cp{2}$ such that there exists a diffeomorphism $\varphi:\cp{2}\to \cp{2}$ identifying 
$P$ with the standard $\rp{2}\subset \cp{2}$.  Moreover, it is such that each holomorphic disk of the 
family represents a generator of $\ho_{2}(\cp{2},P)\cong \bbZ$.  For instance, the standard round metric on $\bbS^{2}$ 
corresponds to a family of holomorphic disks with boundary lying on the standard $\rp{2}\subset \cp{2}$.

In the case of the standard $\rp{2}\subset \cp{2}$, one can check rather directly that there is 
a unique family of holomorphic disks associated to it.  This is because considering the involution
\[
\begin{array}{lccc}
    \rho: & \cp{2} & \to & \cp{2}  \\
          & [x_{0}:x_{1}:x_{2}] & \mapsto & [ \overline{x}_{0}:\overline{x}_{1}:\overline{x}_{2}]
\end{array}
\]
having the standard $\rp{2}$ as a fixed set, we can double any holomorphic disk $D\subset \cp{2}$
with $\pa D\subset \rp{2}$ and $D\setminus \pa D\subset \cp{2}\setminus\rp{2}$ to obtain a holomorphic 
curve
\[
                  \Sigma= D\cup \rho(D)\subset \cp{2}.
\]
In the case of a Zoll family, $D$ has to represent a generator of $\ho_{2}(\cp{2},\rp{2})$, which
means $\Sigma\cong \cp{1}$ is a curve of degree $1$ in $\cp{2}$.  Thus, $\Sigma$ is the zero locus
of a homogeneous polynomial of degree 1.  Furthermore, the coefficients of this polynomial can be 
chosen to be real since $\rho(\Sigma)=\Sigma$.  This indicates that $D$ has to be a member of the 
family of holomorphic disks corresponding to the standard Zoll metric on $\bbS^{2}$.  

However, in the more general situation where $P\subset \cp{2}$ does not come from the fixed set
of an antiholomorphic involution, the above argument has no obvious generalization. In principle, one could therefore
imagine situations where two distincts Zoll metrics on $\bbS^{2}$ would lead to two distinct families
of holomorphic disks with boundary on the same totally real submanifold $P\subset \cp{2}$.  This would
mean the totally real submanifold $P\subset \cp{2}$ would not be sufficient to recover the
Zoll metric from which it comes from.

In this paper, we rule out this possibility by showing that for  a totally real submanifold $P\subset \cp{2}$ such that there exists a diffeomorphism $\varphi:\cp{2}\to \cp{2}$ identifying
$P$ with the standard $\rp{2}$, there is at most one family of holomorphic disks coming from a Zoll metric
on $\bbS^{2}$.  In particular, the twistor correspondence of LeBrun and Mason, which to a Zoll metric
on $\bbS^{2}$ associate a totally real submanifold $P\subset \cp{2}$, is injective.  In fact, we show
that given such a family, any holomorphic disk $D\subset \cp{2}$ with 
$\pa D\subset P$, $D\setminus \pa D\subset \cp{2}\setminus P$ and such that it is a generator of 
$\ho_{2}(\cp{2},P)$ necessarily has to be a member of the family (see theorem~\ref{uzs.1} in \S\ref{ucf.0}).

Our proof is by contradiction.  We suppose that there is such a holomorphic disk $D\subset \cp{2}$
which is not in the family.  Then we show that there exists  a holomorphic disk $D_{b}$ in the family such that
$D$ and $D_{b}$ intersect transversally and have a non-empty intersection in $\cp{2}\setminus P$.  Since
$D$ and $D_{b}$ have boundaries, the intersection number one gets is not a homotopy invariant.  Our 
strategy is to look at the lifts $\widetilde{D}$ and $\widetilde{D}_{b}$ of $D$ and $D_{b}$ in the blow-up
(in the sense of Melrose) 
\[
           \tcp{2}:= [\cp{2};P].
\]
Since the family of holomorphic disks defines a circle fibration on $\pa\tcp{2}$, we have an alternative
way of blowing down $\tcp{2}$.  Let 
\[
       \alpha: \tcp{2}\to Y
\]
denote the corresponding blow-down map.
In this blow-down procedure, each disk becomes a sphere, so eliminating the
boundary.  We can also get a sphere out of $D$ by deforming it before blowing down with $\alpha$.
Doing it carefully, we can insure that the intersection number between $\tD$ and $\tD_{b}$ remains
positive.  One can also check that the corresponding spheres $\alpha(\tD)$ and $\alpha(\tD_{b})$ in the blow-down picture are homologous.  By construction, the oriented intersection number of these two spheres will be positive,
which leads to a contradiction, since the self-intersection number of the corresponding homology class
is shown to be equal to zero.

The paper is organized as follows.  In section~\ref{pre.0}, we review standard results of differential 
topology.  In section~\ref{cf.0}, we introduce the notion of complete family of $J$-holomorphic curves, while in section~\ref{zm.0}, we explain how it naturally arises in the twistor correspondence of 
LeBrun and Mason.  In section~\ref{ucf.0}, we state and prove the main result. 

\section*{Acknowledgement}
The author would like to warmly thank Claude LeBrun for suggesting the problem studied in this paper, for 
many helpful discussions and for giving very nice lectures on twistor spaces.  

\section{Preliminaries}\label{pre.0}

Unless otherwise stated, the manifolds and maps considered in this paper will be assumed to be smooth.
To clarify our conventions with orientation, let us review briefly basic concepts of differential 
topology.  Let $X$ and $Y$ be smooth manifolds, possibly with boundary.  Let $W\subset X$ be a smooth
submanifold, possibly also with a boundary.  Then a smooth map $f: Y\to X$ is said to be 
\textbf{transversal} to the submanifold $W$ if for each $y\in f^{-1}(W)$,
\[
                   f_{*}(T_{y}Y)+ T_{x}W= T_{x}X, \quad x=f(y).
\]
In that case, and when $Y$ and $W$ do not have a boundary, the preimage $S:=f^{-1}(W)$ is
a submanifold of $Y$ of codimension equal to the codimension of $W$ in $X$.  Moreover, if $X,Y$ and
$W$ are oriented manifolds, then $S$ has a natural induced orientation.  More precisely, we will
follow the convention of \cite{Guillemin-Pollack}, p.101. Let $\nb S$ be the normal bundle of $S$ in
$Y$.  Choosing a Riemannian metric on $Y$, we can think of $\nb S$ as a subbundle of 
$\left. TY\right|_{S}$.  Then, for $y\in S$, we have
\begin{gather}
  f_{*}(\nb_{y}S)\oplus T_{x}W= T_{x}X, \quad x=f(y), \label{pre.1}\\
  \nb_{y}S\oplus T_{y}S= T_{y}Y.
\label{pre.2}\end{gather}
Let $\beta_{X}$, $\beta_{Y}$ and $\beta_{W}$ be oriented basis of $T_{x}X$, $T_{y}Y$ and $T_{x}W$
respectively.  The first equation \eqref{pre.1} determines an orientation on $N_{y}S$ by declaring 
a basis $\beta_{N}$ of $\nb_{y}S$ oriented if $(f_{*}(\beta_{N}),\beta_{W})$ is an oriented basis of
$T_{x}X$.  Using the second equation \eqref{pre.2}, this in turn determines an orientation on
$T_{y}S$ by declaring a basis $\beta_{S}$ of $T_{y}S$ oriented if $(\beta_{N},\beta_{S})$ is an
oriented basis of $T_{y}Y$ whenever $\beta_{N}$ is an oriented basis of $\nb_{y} S$.  

A particularly interesting case is when the dimension of $Y$ is the same as the codimension of $W$ in
$X$.  In that case, the submanifold $S:= f^{-1}(W)$ is a discrete set of points which is finite
if $Y$ is compact.  Each of these point has an orientation $\pm 1$ provided by the preimage 
orientation.  We define the \textbf{oriented intersection number} between $f$ and $W$, denoted
$I(f,W)$, to be the sum of their orientations numbers.  Defined this way, it turns out (see for instance
\cite{Guillemin-Pollack}) that the oriented intersection number only depends on the homotopy class
of the map $f$ when both $Y$ and $W$ are compact manifolds without boundary.  More generally, when
$X$, $Y$ and $W$ are not oriented, the intersection number modulo 2 still only depends on the homotopy
class of $f$.  

The transversality of a map is a natural and generic property which is usually easy to satisfy modulo
small perturbations of the map.  In this paper, we will in particular use the following standard
results of differential topology.

\begin{transversality}
Suppose $\phi: Y\to B$ is a (locally trivial) smooth fibration of manifolds with boundary, the 
base $B$ being a compact manifold without boundary.  Let $F:Y\to X$ be a smooth map of manifolds where
$X$ has no boundary, and let $W\subset X$ be a boundaryless submanifold.  If both $F$ and $\pa F$ are
transversal to $W$, then for almost all $b\in B$, both $f_{b}: \left. F\right|_{\phi^{-1}(b)}$ and
$\pa f_{b}:= \left.\pa F\right|_{\pa\phi^{-1}(b)}$ are transversal to $W$.     
\label{pre.3}\end{transversality}  
\begin{remark}
In \cite{Guillemin-Pollack}, the theorem is stated for trivial fibrations, but the proof 
generalizes immediately to locally trivial fibrations.
\label{pre.4}\end{remark}

\begin{extension}
Suppose that $W$ is a closed submanifold of $X$, both boundaryless.  Let $Y$ be a compact 
manifold with boundary and $C\subset Y$ a closed subset.  Let $f:Y\to X$ be a smooth map transversal to $W$ on $C$ such that 
$\pa f: \pa Y\to X$ is transversal to $W$ on $C\cap \pa Y$.  Then there exists a smooth map
$g: Y\to X$ homotopic to $f$ and arbitrarily $\CI$-close to $f$ such that $g$ and $\pa g$ are
transversal to $W$ and $g=f$ in a neighborhood of $C$.
\label{pre.5}\end{extension} 

When $C=\emptyset$, this gives the transversality homotopy theorem.  Another important special case
of the extension theorem is the following.
\begin{corollary}
If for $f: Y\to X$, the boundary $\pa f:\pa Y\to X$ is transversal to $W\subset Y$, then
there exists a map $g: X\to Y$ homotopic to $f$ and arbitrarily $\CI$-close to $f$ such that
$\pa g=\pa f$ and $g$ is transversal to $W$.
\label{pre.6}\end{corollary}
If $X$ is an oriented manifold with boundary, then its boundary has a natural induced orientation
modulo a choice of convention.  This has to do with the fact the normal bundle of $\pa X$ is 
orientable.  We will follow the convention of \cite{Guillemin-Pollack}.  Let $x\in \pa X$ be 
given and let $n_{x}\in T_{x}X$ be a vector pointing outside $X$.  Then we declare a basis
$\beta_{\pa X}$ of $T_{x}\pa X$ oriented if $(n_{x}, \beta_{\pa X})$ is an oriented basis 
of $T_{x}X$.  This defines an orientation on $\pa X$.

In this paper, manifolds with boundary will mostly arise from a blow-up procedure introduced by
Melrose (see for instance \cite{MelroseMWC}).  
Let $Z$ be a closed manifold and  assume $P\subset Z$ is an embedded 
closed submanifold.  Both $Z$ and $P$ are possibly not orientable.  Let 
\begin{equation}
  \nb P := T_{P}Z / TP, \quad \unb P:= (\nb P \setminus \{0_{\nb P}\})/\bbR^{+}
\label{it.1}\end{equation}
be the normal bundle of $P$ in $Z$ with its associated (abstract) unit normal bundle, $0_{\nb P}$ denoting
the zero section of $\nb P$.  Let 
\begin{equation}
    \nu: \nb P \tilde{\longrightarrow} \cU \subset Z
\label{it.2}\end{equation}
be a tubular neighborhood of $P$ in $Z$ identifying the zero section of $\nb P$ with $P\subset Z$.  Let 
$L\to \unb P$ be the (trivial) tautological line bundle, that is, $L\subset \pi^{*} \nb P$, where 
$\pi:\unb\to P$ is the canonical projection, has fibre at $[v]\in \unb P$ given by $\spa(v)$. 
With respect to the canonical `outward pointing' orientation of $L\to \unb P$, let $L^{+}\to \unb P$ be
the $[0,\infty)$-bundle consisting of vectors which are not inward pointing.  Then the map $\nu$ has a 
natural lift $\tilde{\nu}$ so that one has the commuting diagram
\begin{equation}
\xymatrix{ L^{+}\ar[d]^{\psi}\ar[rd]^{\tilde{\nu}} & \\
 \nb P\ar[r]^{\nu} & \cU }
\label{it.3}\end{equation}
where $\psi: L^{+}\to \nb P$ is the natural map defined by
\begin{equation}
    \psi(p,[v],w)= (p,\pi_{*}w), \quad p\in P, \; [v]\in \unb_{p} P, \; 
                                               w\in L^{+}_{(p,[v])}\subset \pi^{*}\nb_{p}P.
\label{it.4}\end{equation}
One can associate a manifold with boundary to the pair $(Z,P)$ by blowing up $Z$ at $P$ in the sense
of Melrose \cite{MelroseMWC}.  Essentially, this amounts to the introduction of polar coordinates.
\begin{definition} The 
\textbf{Melrose's blow up of $Z$} at $P$,
\begin{equation}
  [Z;P] := L^{+} \cup_{\tilde{\nu}} (Z\setminus P)
\label{it.5}\end{equation}
is obtained by gluing $L^{+}$ and $Z\setminus P$ via $\tilde{\nu}$.  It comes together with a blow-down map 
\begin{equation}
              \beta: [Z;P] \to Z    
\label{it.8}\end{equation}
which on $Z\setminus P \subset [Z;P]$ corresponds to the identity map $Z\setminus P \to Z\setminus P$ and
on the zero section of $L^{+}$ corresponds to the canonical projection $\pi: \unb P\to P$.
\label{pre.6}\end{definition}
  In particular, this means the boundary of $[Z;P]$ is canonically identified with $\unb P$,
\begin{equation}
           \pa [Z;P] \cong \unb P.
\label{it.6}\end{equation}
The diffeomorphism class of $[Z;P]$ does not
depend on the choice of the tubular neighborhood map $\nu$ (see \cite{MelroseMWC} \S 5.3).
To lighten the notation, we will sometime write $\tZ:=[Z;P]$ when it is clear what is the
submanifold $P$ where the blow-up is taken.

\section{Complete families of $J$-holomorphic curves}\label{cf.0}

Let $Z$ be a smooth closed 4-manifold  equipped with an almost complex structure $J$.
Suppose that $P\subset Z$ is a totally real submanifold, that is
\begin{equation}
    T_{p}P \cap J(T_{p}P)=0 \quad \forall p\in P, \quad
                                 \left. TZ\right|_{P}= TP \oplus J(TP). 
\label{cf.1}\end{equation}
Let 
\begin{equation}
\xymatrix{ \tSigma \ar@{-}[r]  &  M\ar[d]^{\phi} \\
  &  B }
\label{cf.2}\end{equation}
be a (locally trivial) $\CI$-fibration of manifolds, where the base $B$ is a 
closed manifold, the total space $M$ is a compact manifold with boundary and
the typical fibre $\tSigma$ is a compact Riemann surface with boundary.  
For $b\in B$, let 
\[
           \tSigma_{b}:= \phi^{-1}(b)
\]  
denote the fibre above $b$.  Assume $j\in\CI(M; \End(T(M/B))$ is a 
family of complex structures on the fibres and denote by
\[
             j_{b}\in \CI(\tSigma_{b};\End(T\tSigma_{b}))
\]
the associate complex structure on $\tSigma_{b}$.

\begin{definition}
A \textbf{complete family} of $J$-holomorphic curves with boundary associated to the pair $(Z,P)$ is a 
locally trivial fibration \eqref{cf.2} with total space 
\[
                   M:=\tZ= [Z;P]
\]
such that for each $b\in B$ the restriction of the blow-down map $\beta:\tZ\to Z$ gives
an embedded $J$-holomorphic curve
\[
   \beta_{b}: (\tSigma_{b},j_{b}) \hookrightarrow (Z,J).
\]

\label{cf.11}\end{definition}

When the base $B$ of a complete family \eqref{cf.2} is connected, each $J$-holomorphic curve in
the family defines the same relative homology class in $\ho_{2}(Z,P)$.
In that case, we denote by
$\hcf\in [(\tSigma_{b},\pa\tSigma_{b}); (Z,P)]$ the homotopy class defined by the embedding
$\nu_{b}: \tSigma_{b}\hookrightarrow Z$ for $b\in B$ which maps $\pa \tSigma_{b}$ into
$P$ and we denote by $[\hcf]\in \ho_{2}(Z,P)$ the corresponding relative homology class.  Since
$B$ is connected, the element $\hcf$ does not depend on the choice of $b\in B$.

Given a complete family of $J$-holomorphic curves with \textbf{connected} boundary associated to a 
pair $(Z,P)$, there is an alternative way of blowing down $\tZ:= [Z;P]$ in the sense of Melrose.  Indeed,
the fibration $\phi: \tZ\to B$ associated to such a family gives the boundary $\pa\tZ$ of $\tZ$ the
structure of a circle bundle
\begin{equation}
\xymatrix{ \pa\tSigma\ar@{-}[r] & \pa \tZ \ar[d]^{\pa\phi} \\
          & B }
\label{bd.1}\end{equation}
In other words, this is a $\SO(2)$-principal bundle, so one can associate to it an oriented vector
bundle $V\to B$ of rank two equipped with an inner product on each fibre.  Then the circle
bundle \eqref{bd.1} corresponds to the unit circle bundle of $V$.  The zero section of $V$ gives
an inclusion $B\subset V$.  Consider the blow up \`a la Melrose of $B$ in $V$,
\begin{equation}
    \tV:= [V;B].
\label{bd.2}\end{equation} 
Then the boundary $\pa\tV$ of $\tV$ has an induced circle bundle structure 
\[
\xymatrix{ \bbS^{1} \ar@{-}[r] & \pa \tV \ar[d] \\
          & B }
\]
which is canonically isomorphic to the circle bundle \eqref{bd.1}.  We can choose a 
collar neighborhood
\begin{equation}
   c: \pa\tZ\times [0,1) \hookrightarrow \tZ
\label{bd.3}\end{equation}
of $\pa\tZ$ in $\tZ$ inducing a commutative diagram
\begin{equation}
  \xymatrix{ \tV \ar[r]^{c}\ar[rd] &  \tZ \ar[d]^{\phi} \\
          & B }
\label{bd.4}\end{equation}
with $\pa\tV= \pa \tZ$ and $\tV \setminus \pa \tV\subset \tZ\setminus \pa \tZ$.
Thus, we can consider the compact $\CI$-manifold without boundary 
\begin{equation}
   Y:= V \cup_{c} (\tZ\setminus \pa\tZ)
\label{bd.5}\end{equation}
obtained by identifying $V\setminus B= \tV\setminus \pa\tV$ with 
$c(\tV\setminus \pa \tV)\subset \tZ\setminus \tZ$.  This construction
gives an inclusion $B\subset V\subset Y$ and a canonical identification
\[
          \tZ= [Y; B]        
\]
with blow-down map $\alpha: \tZ\to Y$.  The commutative diagram \eqref{bd.4} induces
a fibration $\psi: Y\to B$ and a commutative diagram
\begin{equation}
\xymatrix{ \tZ \ar[r]^{\alpha}\ar[rd]^{\phi} &  Y \ar[d]^{\psi} \\
          & B }
\label{bd.6}\end{equation}
Under the blow-down map $\alpha$, the Riemann surface with boundary $\tSigma_{b}= \phi^{-1}(b)$
is mapped to its blow-down version $\Sigma_{b}:=\psi^{-1}(b)$, which is a 
compact Riemann surface with empty boundary,
\[
         \tSigma_{b}= [\Sigma_{b};b]
\]
with blow-down map $\alpha_{b}: \tSigma_{b}\to \Sigma_{b}$ where $b\in B$ is seen as a point
on $\Sigma_{b}$ via the inclusion $B\subset V\subset Y$.  
\begin{remark}
When the base $B$ is oriented and when
the fibres of the fibration \eqref{cf.2} are consistently oriented, notice that this induces a natural orientation on $Y$.
\label{orientation}\end{remark}
To summarize, given a complete family of
$J$-holomorphic curves with connected boundary associated to the pair $(Z,P)$, we get a commutative
diagram
\begin{equation}
\xymatrix{ Y\ar[rd]^{\psi}& \tZ\ar[l]_{\alpha}\ar[r]^{\beta}\ar[d]^{\phi} & Z  \\
            & B & Z\setminus P\ar[l]\ar[u] }
\label{bd.7}\end{equation} 
where $\alpha$ and $\beta$ are blow-down maps.  When we have instead a family of $J$-holomorphic curves
with disconnected boundary, there is a similar construction.  We leave the details to the interested reader.

\section{Zoll manifolds and complete families of holomorphic Disks}\label{zm.0}

A Zoll metric on a smooth manifold $M$ is a Riemannian metric $g$ whose geodesics are all simple closed
curves of equal length.  This terminology is in honor of Otto Zoll's result \cite{Zoll} that $\bbS^{2}$
admits many such metrics beside the obvious metrics of constant curvature.  Later on, it was speculated by Funk \cite{Funk} via a perturbation argument and proved
by Guillemin \cite{Guillemin_Zoll} using a type of Nash-Moser implicit function theorem that the general Zoll metric on $\bbS^{2}$ depends modulo isometries and rescaling on one odd function
$f:\bbS^{2}\to \bbR$ (\ie a function such that $f(-\vec{r})=-f(\vec{r})$ for
all $\vec{r}\in \bbS^{2}\subset\bbR^{3}$).  In other words, odd functions on $\bbS^{2}$ describe the tangent bundle of the 
moduli space of Zoll metrics on $\bbS^{2}$ at the standard metric.  In particular, this moduli space is infinite dimensional.  In contrast, on $\rp{2}$, there is only one such metric modulo rescaling and isometries, as was conjectured by Blaschke and proved by Leon Green \cite{Green}.  

More recently, LeBrun and Mason \cite{Lebrun-Mason} were able to recover all these results and to generalize them to Zoll symmetric affine connections (not necessarily the Levi-Civita connection of a Riemannian metric) by using an approach involving twistor theory.  For instance, on $\bbS^{2}$, they established a 
twistor correspondence between Zoll connections and families of holomorphic disks in $\cp{2}$ having
their boundary contained in a totally real submanifold $P\subset \cp{2}$ related by a diffeomorphism of
$\cp{2}$ to the standard $\rp{2}\subset\cp{2}$.

More precisely, the twistor correspondence goes as follows.  Recall first that two torsion-free
affine connections $\nabla$ and $\hat{\nabla}$ on a manifold $M$ are said to be \textbf{projectively equivalent} if they have the same geodesics, considered as unparametrized curves.  A 
\textbf{projective structure} on a smooth manifold is the projective equivalence class $[\nabla]$ of
some affine torsion-free connection.  A projective structure $[\nabla]$ is said to be 
a \textbf{Zoll projective structure} if the image of any maximal geodesic of $\nabla$ is an embedded 
circle $\bbS^{1}\subset M$.

Now, let $[\nabla]$ be a Zoll projective structure on $\bbS^{2}$.  Let $T_{\bbC}\bbS^{2}:= 
T\bbS^{2}\otimes_{\bbR}\bbC$ be the complexified tangent bundle of $\bbS^{2}$.  
Denote by $\bbP T_{\bbC}\bbS^{2}$ the complex projectivization of the complexified tangent
bundle $T_{\bbC}\bbS^{2}$ and by $\bbP T\bbS^{2}$ the real projectivization of the tangent
bundle $T\bbS^{2}$.
Then
\begin{equation}
     \cU= \cZ\setminus Z= \bbP T_{\bbC}\bbS^{2}\setminus \bbP T\bbS^{2}
\label{vvv.1}\end{equation}
can be identified with the space of all pointwise almost-complex structures on $\bbS^{2}$.  Indeed,
given a pointwise complex structure $J$ on $T_{s}\bbS^{2}$, $s\in \bbS^{2}$, we let
$[(v+ iJv)]\in \cU$ be the corresponding element in $\cU$ where $v\in T_{s}\bbS^{2}$ is any
non-zero vector.  Thus,
\[
             \cU= \cU_{+}\cup \cU_{-}
\]
where $\cU_{+}$ (respectively $\cU_{-}$) consists of the point-wise almost complex structures
compatible (respectively incompatible) with the orientation of $\bbS^{2}$.  Both 
$\cU_{+}$ and $\cU_{-}$ are connected.  Consider the manifold with boundary
\[
       \cZ_{+}:= \cU_{+}\cup \bbP T\bbS^{2}, \quad \pa \cZ_{+} = \bbP T\bbS^{2}.
\] 
The bundle projection $\bbP T_{\bbC}\bbS^{2} \to \bbS^{2}$ induces by restriction a 
locally trivial fibration of disks on $\cZ_{+}$,
\begin{equation}
\xymatrix{  \bbD\ar@{-}[r] & \cZ_{+}\ar[d]^{\mu} \\
           & \bbS^{2}  }
\label{pre.7}\end{equation}
with zero section having self intersection 4.
Now, the geodesics of the Zoll structure naturally lift to $\bbP T\bbS^{2}$ to 
give a foliation of $\bbP T\bbS^{2}$ by circles.  As discussed in \cite{Lebrun-Mason}, this
is in fact a locally trivial fibration
\begin{equation}
\xymatrix{  \bbS^{1} \ar@{-}[r] & \pa\cZ_{+}\ar[d]^{\nu} \\
           & \rp{2}  }
\label{pre.77}\end{equation}
over $\rp{2}$ which is isomorphic to the unit tangent bundle 
$\bbS T\rp{2}\to \rp{2}$.  Thus, as in the discussion following \eqref{bd.1}, there is 
a corresponding blow-down procedure (in the sense of Melrose)
\begin{equation}
  \beta: \cZ_{+}\to Z
\label{pre.8}\end{equation}
where $Z$ is a compact manifold with no boundary.  It turns out to be diffeomorphic to $\cp{2}$.  In
fact, $Z$ comes equipped with a natural complex structure induced by the Zoll projective structure
$[\nabla]$ on $\bbS^{2}$.  The complex structure arises as follows.  The $(0,1)$-tangent bundle of the 
fibres of the fibration \eqref{pre.7} defines a complex line bundle
\begin{equation}
   L_{1}:= T^{0,1}(\cZ_{+}/\bbS^{2})\subset T_{\bbC}\cZ_{+}= T\cZ_{+}\otimes_{\bbR}\bbC.
\label{dist.1}\end{equation}  
On the other hand, a choice of representative $\nabla$ for the Zoll projective structure $[\nabla]$
defines a horizontal lift of $T_{\bbC}\bbS^{2}$ to $\bbP T_{C}\bbS^{2}$, and so a horizontal lift
\[
       H_{\bbC}\cong \mu^{*}T_{\bbC}\bbS^{2}
\]
to $\cZ_{+}$.  In particular, one can define a tautological line subbundle $L_{2}\subset H_{\bbC}$ by
\[
    \left.L_{2}\right|_{[w]}= (\mu_{* [w]})^{-1}(\Span w)
\]
for $[w]\in \cZ_{+}$.  In this way, one gets a $\CI$-distribution 
\begin{equation}
      \mathbf{D}:= L_{1}\oplus L_{2} \subset T_{\bbC} \cZ_{+}.
\label{dist.2}\end{equation}
LeBrun and Mason showed \cite{Lebrun-Mason} that the distribution $\mathbf{D}$
is involutive on $\cZ_{+}$, in the sense that
\[
     [\mathcal{C}^{1}(\mathbf{D}),\mathcal{C}^{1}(\mathbf{D})]\subset \mathcal{C}^{0}(\mathbf{D}),
\]
and that it is independent of the choice of representative $\nabla$ of the Zoll projective structure
$[\nabla]$.  Since
\[
    \dim( \mathbf{D}_{z}\cap\overline{\mathbf{D}}_{z})= \left\{
     \begin{array}{ll}
       0, & z\in \cZ_{+}\setminus\pa\cZ_{+}, \\
       1, & z\in \pa\cZ_{+},
     \end{array} \right.
\]
one can apply the Newlander-Nirenberg theorem to get an integrable 
complex structure away from the boundary.  
Notice also that definition \eqref{dist.2} works also perfectly well if the projective
structure is not Zoll.  However, when it is Zoll, one can consider the image of $\mathbf{D}$ under the blow-down
map \eqref{pre.8}.  LeBrun and Mason showed that 
\[
    \beta_{*}(\mathbf{D})\subset T_{\bbC}Z
\]
is still involutive and such that $\beta_{*}(\mathbf{D})\cap\beta_{*}(\overline{\mathbf{D}})=0$, so that it can be interpreted as the $(0,1)$-tangent bundle of an appropriate
(integrable) complex structure.

With this complex structure, $Z$ is biholomorphic to $\cp{2}$.  Most
importantly,  the disk fibration \eqref{pre.7} together with the blow-down map \eqref{pre.8} define
a complete family of holomorphic disks associated to the pair $(\cp{2},P)$ where 
$P:=\beta(\pa \cZ_{+})\subset \cp{2}$ is diffeomorphic to $\rp{2}$.
With this notation, the manifold with boundary $\cZ_{+}$ can be naturally identified with the
blow-up in the sense of Melrose of $\cp{2}$ at $P$,
\[
                 \tcp{2}:= [\cp{2};P].
\]  
In this context, the commutative diagram \eqref{bd.7} takes the form
\begin{equation}
\xymatrix{ Y\ar[rd]^{\psi}& \tcp{2}\ar[l]_{\alpha}\ar[r]^{\beta}\ar[d]^{\mu} & \cp{2}  \\
            & \bbS^{2} & \cp{2}\setminus P\ar[l]\ar[u] }
\label{bd.7b}\end{equation}
where $Y$ is the $\bbS^{2}$-bundle associated to the disk fibration \eqref{pre.7}.

In this correspondence, holomorphic disks correspond to points on $\bbS^{2}$, while points on $P$ correspond
to geodesics of the projective structure of $[\nabla]$ on $\bbS^{2}$.  Moreover, the points of intersections
between two holomorphic disks corresponds to the geodesics passing through the two points on $\bbS^{2}$ associated to the two holomorphic disks.
In \cite{Lebrun-Mason}, it is shown
that there exists a diffeomorphism $\varphi: \cp{2}\to \cp{2}$ mapping $P$ diffeomorphically onto
the standard $\rp{2}\subset \cp{2}$.  The complete family arising in this way is also such that
the relative homology class $[h_{\mu}]$ represented by a disk of the family is a generator of 
$\ho_{2}(\cp{2},P)\cong \bbZ$.  

For the standard round metric on $\bbS^{2}$, the totally real submanifold $P\subset \cp{2}$ coming from the twistor correspondence is the standard $\rp{2}\subset \cp{2}$.  It is precisely
the set of fixed points of the involution
\begin{equation}
  \begin{array}{lccc}
    \rho: & \cp{2} & \to & \cp{2}  \\
          & [x_{0}:x_{1}:x_{2}] & \mapsto & [ \overline{x}_{0}:\overline{x}_{1}:\overline{x}_{2}]
\end{array}
\label{invol.1}\end{equation}
Each projective line $\Sigma\subset \cp{2}$ invariant under this involution is such that
$\Sigma\cap \rp{2}$ is a a real projective line in $\rp{2}$.  In particular, $\Sigma= D_{+}\cup
D_{-}$ is the union of two holomorphic disks $D_{+}$ and $D_{-}$ with common boundary
\[
    \pa D_{+}= \pa D_{-} = \Sigma\cap \rp{2}
\]
contained in $\rp{2}$.  The holomorphic disks of the complete family associated to the round
metric are precisely those arising in this way.  Thus, the holomorphic disks of this 
complete family come into pairs $D_{b}, \rho(D_{b})$.  This extra symmetry is very specific
to the standard round metric and does not arise for a general Zoll metric.  It corresponds to the fact that for each point $\vec{r}$ on $\bbS^{2}$,  the points $\vec{r}$ and $-\vec{r}$ 
are conjugate and have exactly the same set of geodesics
passing through them (the geodesics being great circles in this case).

The following lemma will turn out to be useful in the next section.
\begin{lemma}
Let $[\nabla]$ be a Zoll projective structure on $\bbS^{2}$ and let $p\in \bbS^{2}$ be given.  Then,
for $q\in \bbS^{2}\setminus\{p\}$ sufficiently close to $p$, there exists a unique geodesic circle
of $[\nabla]$ joining $p$ and $q$.
\label{geod.1}\end{lemma}
\begin{proof}
The diagrams \eqref{pre.7} and \eqref{pre.77} combine to give a diagram
\begin{equation}
\xymatrix{  & \bbP T\bbS^{2}\ar[ld]_{\mu}\ar[rd]^{\nu} &  \\
              \bbS^{2} & & \rp{2}
}
\label{geod.2}\end{equation}
where we have used the identification $\bbP T\bbS^{2}\cong \pa \cZ_{+}$.
The maps $\mu$ and $\nu$ are both locally trivial circle fibrations, $\mu$ being the standard 
bundle projection while the fibres of $\nu$ are the geodesics of $[\nabla]$ lifted to 
$\bbP T\bbS^{2}$.  In $\bbP T\bbS^{2}$, consider the union
\[
       \hat{X}= \nu^{-1}( \nu[ \mu^{-1}(p)])
\]
of the lift of the geodesics of $[\nabla]$ passing through $p$.  As discussed in (\cite{Lebrun-Mason},
lemma 2.8), $\hat{X}$ is a compact smooth $2$-manifold which can be blown down (in the usual projective sense)
along $\mu^{-1}(p)$ since the normal bundle of $\mu^{-1}(p)\cong \rp{1}$  in
$\hat{X}$ is isomorphic to the universal line bundle of $\rp{1}$.  The blow-down of $\hat{X}$ produces
a new manifold $X$ and the map $\mu$ induces a smooth map
\[
         \varrho: X\to \bbS^{2}.
\]
If $\check{p}\in X$ is the blow-down of $\mu^{-1}(p)$, then $\varrho$ is modeled on the exponential
map of $\nabla$ near $\check{p}$.  In particular, there exists an open neighborhood $\cU$ of $p$
in $\bbS^{2}$ such that 
\[
        \varrho: \varrho^{-1}(\cU)\to \cU
\]
is a diffeomorphism.  But then, since $[\nabla]$ is a Zoll projective structure, 
\[
     p \not\in \varrho( X\setminus \varrho^{-1}(\cU)).
\]
Thus, there exists an open neighborhood $\cV$ of $p$ contained in $\cU$ such that 
\[
     \cV\cap \varrho(X\setminus \varrho^{-1}(\cU))= \emptyset.
\]
The result then follows by noticing that for all $q\in \cV\setminus\{p\}$, we have
$\varrho^{-1}(q)\subset \varrho^{-1}(\cU)$, so that for each $q\in \cV\setminus\{p\}$,
there is a unique geodesic circle of $[\nabla]$ joining $p$ to $q$.
\end{proof}

\section{Uniqueness of Zoll families}\label{ucf.0}

In this section, we will establish a uniqueness statement concerning complete families of 
$J$-holomorphic disks corresponding to projective Zoll structures on $\bbS^{2}$ via the 
twistor correspondence of LeBrun and Mason.

\begin{definition}
Let $P\subset \cp{2}$ be a totally real submanifold of $\cp{2}$ with respect to the standard 
complex structure $J$ of $\cp{2}$ and suppose that there exists a diffeomorphism 
$\varphi: \cp{2}\to \cp{2}$ identifying $P$ with the standard $\rp{2}\subset\cp{2}$.  Then a 
complete family $\phi: \tcp{2}\to \bbS^{2}$ associated to the totally real embedding $P\subset \cp{2}$
is called a \textbf{Zoll family of holomorphic disks} if it comes from a projective Zoll structure
on $\bbS^{2}$ via the twistor correspondence of LeBrun and Mason.
\label{zfhd.1}\end{definition}

\begin{theorem}
Let $P\subset \cp{2}$ be a totally real submanifold of $\cp{2}$ with respect to the standard 
complex structure $J$ of $\cp{2}$ and suppose that there exists a diffeomorphism 
$\varphi: \cp{2}\to \cp{2}$ identifying $P$ with the standard $\rp{2}\subset\cp{2}$.  Suppose
that $\phi: \tcp{2}\to \bbS^{2}$ is a Zoll family associated to the pair $(\cp{2},P)$.  In particular $[\hcf]\in \ho_{2}(\cp{2},P)$ is a generator of $\ho_{2}(\cp{2},P)$.  Let  
$D\subset \cp{2}$ be any embedded holomorphic disk such that $\pa D\subset P$, $D\setminus \pa D \subset \cp{2}\setminus P$ with
$[D]\in\ho_{2}(\cp{2},\rp{2})\cong \bbZ$ a generator.  Then $[D]=[\hcf]$ and there exists $b\in \bbS^{2}$
such that $\beta_{b}: \phi^{-1}(b)\to D$ is a biholomorphism where $\beta:\tcp{2}\to \cp{2}$ is the
blow-down map and $\beta_{b}:= \left. \beta \right|_{\phi^{-1}(b)}$.
\label{uzs.1}\end{theorem}
An immediate consequence of this theorem is the following.
\begin{corollary}
Let $P\subset \cp{2}$ be a totally real submanifold of $\cp{2}$ with respect to the standard 
complex structure $J$ of $\cp{2}$ and suppose that there exists a diffeomorphism $\varphi:\cp{2}\to \cp{2}$ identifying $P$ with the standard $\rp{2}\subset \cp{2}$.  Then there is at most 
one Zoll family of holomorphic disks $\phi: \tcp{2}\to \bbS^{2}$ associated to the pair
$(\cp{2},P)$.
\label{uzs.2}\end{corollary}

We will prove this theorem in a series of lemmas.
\begin{lemma}
Suppose that $\phi:\tcp{2}\to \bbS^{2}$ is a Zoll family of holomorphic disks.  Suppose also that $D\subset \cp{2}$ is a holomorphic disk with $\pa D\subset P$, $(D\setminus \pa D)\cap P=\emptyset$.  Then either $D$ is an element of the complete family or else there 
exists $b\in B$ such that $\beta_{b}: \phi^{-1}(b)\to \cp{2}$ is transversal to $D$ and has a 
non-empty intersection with $D\setminus \pa D$.
\label{uzs.3}\end{lemma}
\begin{proof}

Given a holomorphic disk $D\subset \cp{2}$ with $\pa D\subset P$ and $(D\setminus \pa D)\cap P=\emptyset$, consider its lift 
$\tD:= \beta^{-1}(D)$ to $\tcp{2}$.  Since $\pa\beta:\pa\tcp{2}\to P$ is a 
submersion, by the transversality theorem~\ref{pre.3}, for almost all $b\in \bbS^{2}$, 
$\pa\beta:\pa\phi^{-1}(b)\to P$ is transversal to $\pa D\subset P$.  Since $J$
identifies $TP$ with $\nb P$ and since $D\subset \cp{2}$ and $\beta(\phi^{-1}(b))\subset \cp{2}$
are embedded holomorphically, this means that for almost all $b\in \bbS^{2}$, 
$\pa \tD\cap \phi^{-1}(b)=\emptyset$.  Applying the transversality theorem to the identity map
$\tcp{2}\setminus\pa\tcp{2}\to \tcp{2}\setminus\pa\tcp{2}$, we see also that for almost all $b\in B$,
$\phi^{-1}(b)\setminus \pa \phi^{-1}(b)$ is transversal to $\tD\setminus\pa \tD$.  Thus, for almost all
$b\in \bbS^{2}$, $\phi^{-1}(b)$ is transversal to $\tD$.  Again, since $J$ identifies $TP$ with
$\nb P$ and since $D$ and $\beta(\phi^{-1}(b))$ are embedded holomorphically, this also
means that $\beta_{b}: \phi^{-1}(b)\to \cp{2}$ is transversal to $D\subset \cp{2}$ for almost
all $b\in \bbS^{2}$.

Now, suppose that $D\subset \cp{2}$ does not belong to the complete family
$\phi: \tcp{2}\to \bbS^{2}$.  Since $\beta: \tcp{2}\setminus\pa\tcp{2}\to
\cp{2}\setminus P$ is a diffeomorphism, this means that there exists $x\in D\setminus \pa D$
such that $\ker(\phi_{*})_{x}\cap T_{x}D\ne T_{x}D$.  Since $D$ and 
$\beta(\phi^{-1}(b_{x}))\subset\cp{2}$ are  holomorphic disks where $b_{x}=\phi(x)$, this means in fact that $\ker(\phi_{*})_{x}\cap T_{x}D= \{0\}$.  In particular, $\phi_{*}(T_{x}D)= T_{b_{x}}\bbS^{2}$, so
for all $b\in \bbS^{2}$ sufficiently close to $b_{x}$, we will have that 
\[
          \beta(\phi^{-1}(b))\cap (D\setminus \pa D) \ne \emptyset.
\]
In particular, we can choose such a $b$ so that $\beta:\phi^{-1}(b)\to \cp{2}$ is 
transversal to $D$.
\end{proof}

To prove theorem~\ref{uzs.1}, we need to show that the second possibility in lemma~\ref{uzs.3} cannot
occur, namely, that the holomorphic disk $D\subset \cp{2}$ cannot intersect transversely a holomorphic
disk $\beta(\phi^{-1}(b))\subset \cp{2}$ of the complete family $\phi: \tcp{2}\to \bbS^{2}$ with
a non-empty intersection in the interior.  The idea of the proof will be to derive a contradiction  
from the existence of such a disk.

Hence, suppose that there is a holomorphic disk $i:D\hookrightarrow \cp{2}$ such that 
$i(\pa D)\subset P$, $i(D\setminus \pa D)\subset \cp{2}\setminus P$ and such that there exists
a complete family $\phi: \tcp{2}\to \bbS^{2}$ of holomorphic disks as in theorem~\ref{uzs.1}.  Suppose
also that  $[D]\in \ho_{2}(\cp{2},P)$ is a generator, that there exists $b\in \bbS^{2}$ such that $D_{b}:= \beta(\phi^{-1}(b))\subset \cp{2}$ 
intersects transversely with $D\subset\cp{2}$ and that 
$(D\setminus \pa D)\cap (D_{b}\setminus \pa D_{b})\ne \emptyset$.  Then the two disks
$D$ and $D_{b}$ have a well-defined intersection number
\[
    I(D,D_{b})= \# \{ D_{b}\cap D \subset \cp{2} \}.
\]
Since $D$ and $D_{b}$ are $J$-holomorphic curves, each intersection point $x\in D_{b}\cap D$
has a positive induced orientation.
Hence, counting with 
orientation, we get the same
intersection number as in $I(D,D_{b})$.  In differential topology, the purpose of counting 
intersection points taking into account the orientation is to get a homotopy invariant and sometime
even a homological invariant.  However, in our case, due to the presence of boundaries, there is 
no direct hope for such a homotopy invariance.  Indeed, one can easily construct examples where 
an intersection point in the interior is gradually moved to the boundary via a smooth homotopy until
it completely disappears.

To get an intersection number which is a homotopy invariant, our strategy consists in getting rid
of the boundaries of $D$ and $D_{b}$ by using the 
blow-down map $\alpha:\tcp{2}\to Y$ of \eqref{bd.6}.  This requires first to deform $i:D\hookrightarrow
\cp{2}$ to put it in a suitable position with respect to this blow-down map.  

Now, since both $D$ and $D_{b}$ give a generator of $\ho_{2}(\cp{2},P)$ by hypothesis, their
boundaries in $P\cong \rp{2}$ give the generator of $\ho_{1}(P)\cong \bbZ_{2}$.  In particular,
this means that $D$ and $D_{b}$ necessarily intersect on their boundaries (in $\cp{2}$) and that
\[
    I(\pa D, \pa D_{b})= \# \{ \pa D\cap \pa D_{b}\}= 1 \mod 2.
\] 
We would like to deform the embedding $i:D\hookrightarrow \cp{2}$ (not necessarily through holomorphic
embeddings) in such a way that $\pa D$ intersects $\pa D_{b}$ in exactly one point in $P$.

\begin{lemma}
Let $\gamma_{i}: \bbS^{1}\hookrightarrow \rp{2}$, $i\in\{1,2\}$ be two embedded circles in 
$\rp{2}$ which intersect transversely and such that each of them gives a generator of 
$\pi_{1}(\rp{2})\cong \bbZ_{2}$.  Then there exists a smooth isotopy of embeddings
\[
              F: \bbS^{1}\times [0,1] \to \rp{2}
\]
such that
\begin{itemize}
\item[(i)] $F(\cdot,0)= \gamma_{1}$;
\item[(ii)] The function $f(t):= \#\{ F(\cdot,t)\cap \gamma_{2}\}$, $t\in [0,1]$ is
         decreasing on $[0,1]$ and $f(1)=1$;
\item[(iii)] Except for a finite subset of $t\in [0,1]$ not containing $0$ or $1$, 
         $F_{t}:\bbS^{1}\hookrightarrow \rp{2}$ with $F_{t}(x):= F(x,t)$ is transversal to 
         $\gamma_{2}$.                  
\end{itemize}
\label{uzs.4}\end{lemma}
\begin{proof}
Since $\gamma_{1}$ and $\gamma_{2}$ both give the generator of $\pi_{1}(\rp{2})\cong\bbZ_{2}$, we 
know that their intersection number is odd, that is,
\[
         I(\gamma_{1},\gamma_{2}):= \#\{\gamma_{1}\cap \gamma_{2}\}= 2n+1
\]
for some $n\in \bbN\cup \{0\}$.  To construct the isotopy $F$, we can proceed by recurrence on
$n\in \bbN\cup\{0\}$.  If $n=0$, we can take the identity isotopy
\[
       F(x,t)= \gamma_{1}(x), \quad \forall x\in \bbS^{1}, \; t\in [0,1].
\]
If $n\ge 1$, let $p\in \rp{2}$ be one of the intersection points of $\gamma_{1}$ and
$\gamma_{2}$.  Identify $\bbS^{1}$ with the unit circle in $\bbC$,
\[
        \bbS^{1}= \{e^{i\theta}\in \bbC \quad | \quad \theta\in \bbR\}.
\]
Without loss of generality, we can assume that for $\theta=0$, we have
\[
     \gamma_{1}(e^{i0})=\gamma_{1}(1)=p= \gamma_{2}(1).
\]
Let $\theta_{1}\in (0,2\pi)$ be the smallest angle such that 
\[
     \gamma_{1}(e^{i\theta_{1}})\cap \gamma_{2}(\bbS^{1})\ne \emptyset.
\]
This angle exists because $n\ge 1$.  Let $q= \gamma_{1}(e^{i\theta_{1}})$ be the
corresponding point of intersection.  Then there exists $\theta_{2}\in (0,2\pi)$ such that
$\gamma_{2}(e^{i\theta_{2}})=q$.  Consider the two loops
\begin{equation*}
\begin{gathered}
\nu_{+}= \left\{  \begin{array}{ll}
                   \gamma_{1}\circ\exp( 2it\theta_{1}), & 0\le t\le \frac{1}{2},  \\
                   \gamma_{2}\circ\exp(i\theta_{2}+i(2t-1)(2\pi-\theta_{2})), & \frac{1}{2}\le t \le 1,
                  \end{array}  \right.  \\
               \nu_{-}= \left\{  \begin{array}{ll}
                   \gamma_{1}\circ\exp( 2it\theta_{1}), & 0\le t\le \frac{1}{2},  \\
                   \gamma_{2}\circ\exp\left(i\theta_{2}(-2t+2)\right), & \frac{1}{2}\le t \le 1.
                  \end{array}  \right.   
\end{gathered}
\end{equation*}
As elements of the fundamental group of $\pi_{1}(\rp{2})\cong \bbZ_{2}$, we have
\[
    [\nu_{+}]\circ[\nu_{-}]^{-1}= [\gamma_{2}]=1 \in \bbZ_{2}.
\]
Thus, one of the loops $\nu_{+}$ and $\nu_{-}$ is contractible in $\rp{2}$, while the other
is the generator of $\pi_{1}(\rp{2})$.  Let us denote by $\nu$ the contractible loop among the
two of them, so
\begin{equation*}
\begin{array}{ll}
   \nu(t)= \gamma_{1}\circ \exp( 2it\theta_{1}), & 0\le t\le \frac{1}{2}, \\
   \nu(t)\in \gamma_{2}(\bbS^{1}),  & \frac{1}{2}\le t\le 1,
\end{array}
\end{equation*}
and $\nu: [0,1]\to \rp{2}$ is embedded in $\rp{2}$.  Since $\nu$ is contractible and embedded, the loop
$\nu$ divides $\rp{2}$ into two distinct regions, the `inside' and the `outside' (see figure~\ref{figure1}).  For instance, one can
see this by lifting $\nu$ to the universal cover $\bbS^{2}$ of $\rp{2}$ and by applying the classical
Jordan curve theorem.  Let $D\subset \rp{2}$ denote one of these regions so that 
$\pa D= \nu([0,1])$, more precisely the one such that $\gamma_{1}(e^{i\theta})\notin D$ 
for $-\epsilon < \theta<0$ and $\epsilon>0$ small.  
If for all $\theta\in (\theta_{1}, 2\pi)$, $\gamma_{1}(e^{i\theta})\notin D$, then for  some $\epsilon>0$,
we can find an isotopy of embeddings from $\gamma_{1}$ to $\tilde{\gamma}_{1}$ which is the 
identity on $e^{i\theta}$ for $\theta\in( \theta_{1}+\epsilon, 2\pi -\epsilon)$ with 
$\epsilon>0$ small and such that $\tilde{\gamma_{1}}$ is transversal to $\gamma_{2}$ and
\[
   \tilde{\gamma}_{1}(\bbS^{1})\cap D= \emptyset, \quad 
   \tilde{\gamma}_{1}(\bbS^{1})\cap \gamma_{2}(\bbS^{1})= (\gamma_{1}(\bbS^{1}) \cap \gamma_{2}(\bbS^{1}))
   \setminus \{p,q\}.
\]  
Taking $\tilde{\gamma}_{1}$ instead of $\gamma_{1}$ decreases $n$ by one in this case and so achieves the 
inductive step of our construction.
\begin{figure}[h]
\begin{minipage}[t]{5.5cm}
\resizebox{4cm}{!}{
\includegraphics[0.5in,1in][5in,4in]{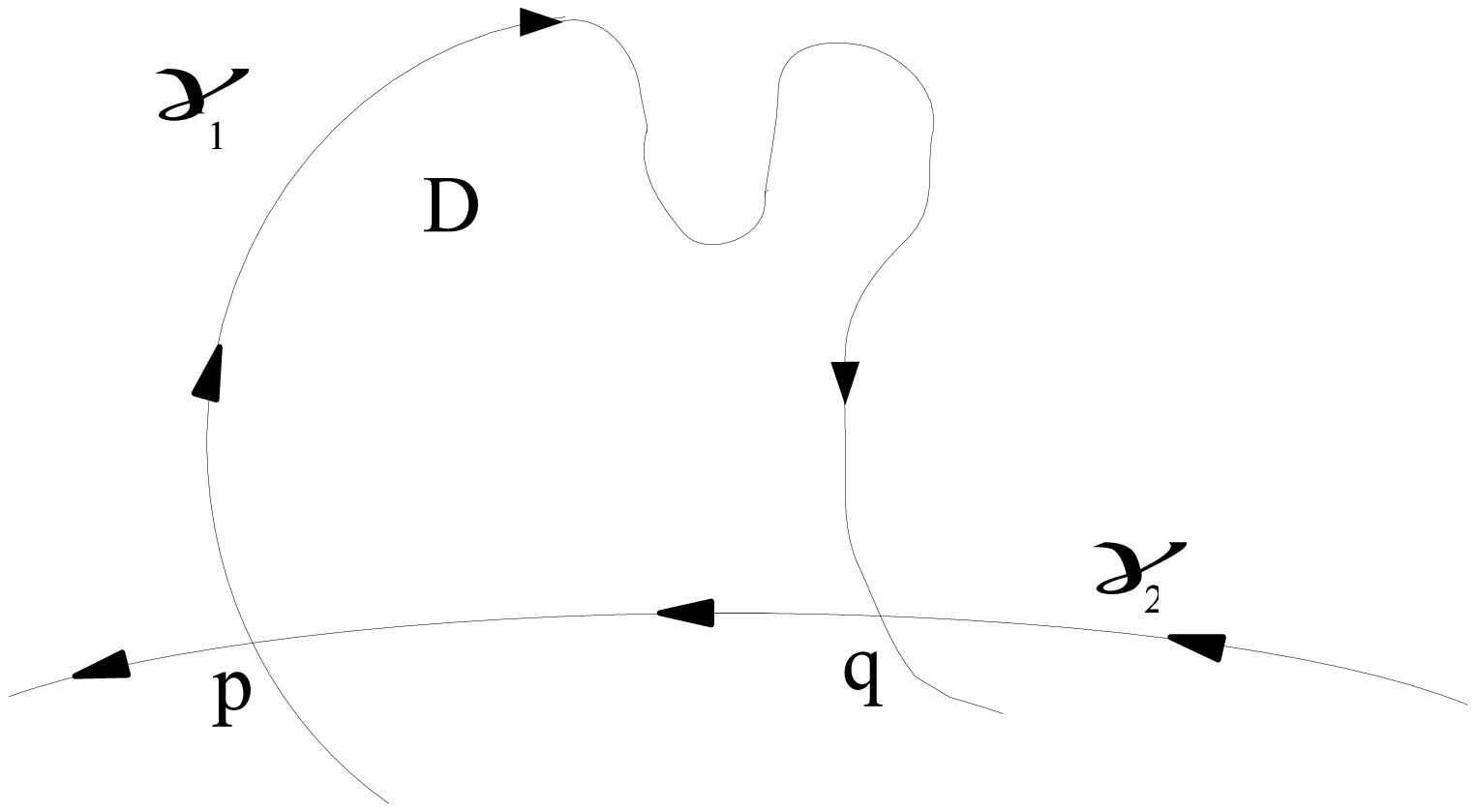}}
	\caption{}
	\label{figure1}
\end{minipage}
\hfill
\begin{minipage}[t]{5.5cm}
\resizebox{4cm}{!}{\includegraphics[0.5in,1in][5in,4in]{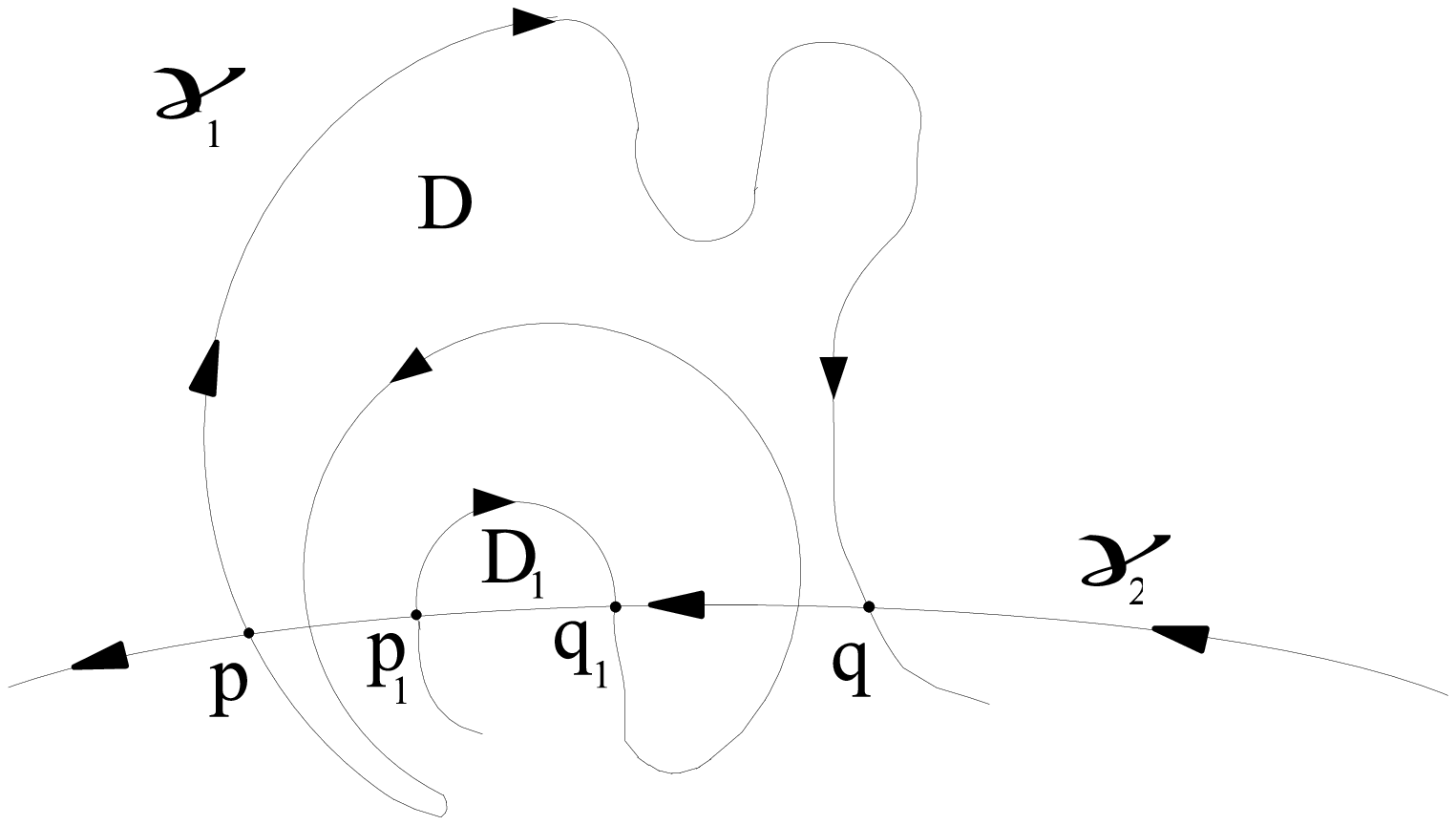}}
	\caption{}
	\label{figure2}
\end{minipage}	
\end{figure}

Otherwise, $\gamma_{1}$ intersects $\nu((\frac{1}{2},1))$ an \textbf{even} number of times
since $\gamma_{1}$ has to get out of $D$ each time it gets inside (see figure~\ref{figure2}).  In particular, we can find
$p_{1}, q_{1}\in \nu(\frac{1}{2},1)$ with
\begin{equation*}
  \begin{array}{lll}
   \gamma_{1}(e^{i\varphi_{1}})= p_{1}, & \gamma_{1}(e^{i\varphi_{2}})=q_{1}, &
                 \theta_{1}< \varphi_{1} < \varphi_{2} < 2\pi,  \\
      \nu(t_{1})=p_{1}, & \nu(t_{2})= q_{1}, &  t_{1},t_{2}\in (\frac{1}{2},1),
  \end{array}
\end{equation*}
and such that the loop
\begin{equation*} \nu_{1}(t):= \left\{
\begin{array}{ll}
  \gamma_{1}\circ\exp( i\varphi_{1}+ 2it(\varphi_{2}-\varphi_{1})), & t\in [0,\frac{1}{2}], \\
  \nu( (2t-1)(t_{1}-t_{2})+ t_{2}) & t\in [\frac{1}{2},1],
\end{array}\right.
\end{equation*}
divides $D$ into two distinct regions $D_{1}$ and $D_{1}'$, $D_{1}$ denoting the region
not containing $\gamma_{1}(e^{i\theta})$ for $\theta\in [0,\theta_{1}]$.  By construction,
\[
    \#\gamma_{1}(\bbS^{1})\cap \nu_{1}((\frac{1}{2},1)) < \# \gamma_{1}(\bbS^{1})\cap
                  \nu((\frac{1}{2},1)).
\]
Since $\# \gamma_{1}(\bbS^{1})\cap \nu((\frac{1}{2},1))$ is finite, we see that proceeding 
recursively, we can in fact assume that 
\[
               \gamma_{1}(\bbS^{1})\cap \nu_{1}((\frac{1}{2},1))= \emptyset.
\]
In that case, we can find an isotopy of embeddings from $\gamma_{1}$ to $\tilde{\gamma}_{1}$ 
which is the identity on $e^{i\theta}$ for $\theta\in (\varphi_{1}-\epsilon, \varphi_{2}+\epsilon)$ and
$\epsilon>0$ small enough, and such that $\tilde{\gamma}_{1}$ is transversal to $\gamma_{2}$ with
\[
   \tilde{\gamma}_{1}(\bbS^{1})\cap D_{1}=\emptyset, \quad 
    \tilde{\gamma}_{1}(\bbS^{1})\cap \gamma_{2}(\bbS^{1})= \gamma_{1}(\bbS^{1}) \cap
     \gamma_{2}(\bbS^{1})\setminus \{p_{1},q_{1}\}.
\]
Replacing $\gamma_{1}$ by $\tilde{\gamma}_{1}$ decreases $n$ by $1$ in this case, also achieving
the inductive step of our construction, which concludes the proof.
\end{proof}

\begin{figure}[h]
\resizebox{4cm}{!}{
\includegraphics[1in,1in][6in,5in]{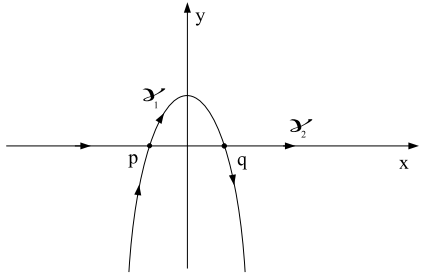}}
	\caption{local model}
	\label{figure3}
\end{figure}

\begin{remark}
Later on, we will need a more precise description of the homotopy achieving the inductive
step of the previous proof.  To this end, we choose a particular \textbf{local model} for the isotopy $F$ around
points where $F_{t}$ is not transversal to $\gamma_{2}$, that is, a local model for each of the inductive steps of the previous proof.  Let $p$ and $q$ (or $p_{1}$ and $q_{1}$) be the two points
involved in the inductive step.  After performing a preliminary isotopy of embeddings to $\gamma_{1}$,
we can assume there are  suitable coordinates $(x,y)$ on an open set $\cU\subset P$ containing $p$ and $q$ such that in these coordinates, $\gamma_{1}$ and $\gamma_{2}$ are locally given by (see figure~\ref{figure3})
\begin{equation}
\begin{gathered}
  \gamma_{1}(\bbS^{1})\cap \cU= \{ u\in \cU \quad | \quad y(u)= -\epsilon x^{2}(u) + \delta\}, \\
  \gamma_{2}(\bbS^{1})\cap \cU= \{ u\in \cU \quad | \quad y(u)=0 \}, \\
  y(p)=y(q)=0, \quad x(p)= -\sqrt{\frac{\delta}{\epsilon}}, \quad x(q)= \sqrt{\frac{\delta}{\epsilon}}.
\end{gathered}
\label{lm.2}\end{equation}
Assume also that the open set $\cU$ contains the disk of radius $\delta+ \sqrt{\frac{\delta}{\epsilon}}$ centered at $(x,y)=(0,0)$.  
To achieve the inductive step, we can then consider the isotopy given locally by translating in the $y$ direction
\begin{equation}
\begin{aligned}
   x(F(e^{i\theta},t))&= x(\gamma_{1}(e^{i\theta})), \\
   y(F(e^{i\theta},t))&= y(\gamma_{1}(e^{i\theta}))-t
\end{aligned}
\label{lm.3}\end{equation} 
with $t\in[0,2\delta]$.  Of course, we need also to adjust this isotopy near $\pa \cU$ so that 
it continues to agree with $\gamma_{1}$ outside of $\cU$.  For each time $t$ where $F_{t}$ is not
transversal to $\gamma_{2}$, we will assume the isotopy achieving the inductive step is given
locally by \eqref{lm.3} in a suitable choice of coordinates around the point where $F_{t}$ is not
transversal to $\gamma_{2}$.
\label{lm.1}\end{remark}

Going back to the holomorphic disks $D$ and $D_{b}$ intersecting transversely and having a non-empty
intersection in $\cp{2}\setminus P$, we can apply lemma~\ref{uzs.4}  to the embedded circles
\[
     \pa i: \pa D\hookrightarrow P, \quad \pa \beta_{b}: \pa \tD_{b}\hookrightarrow P
\]
where $\tD_{b}$ is the lift of $D_{b}$ to $\tcp{2}$.
Thus, identifying $\pa D$ with $\bbS^{1}$, there exists a smooth isotopy of embeddings
\[
           \pa I: \bbS^{1}\times [0,1] \to P
\]
satisfying the three properties of lemma~\ref{uzs.4}, namely
\begin{enumerate}
\item[(i)] $\pa I(\cdot,0)= \pa i$;
\item[(ii)] The function $f(t):= \#(\pa I_{t}(\bbS^{1})\cap \pa D_{b})$ is decreasing on $[0,1]$ and
      $f(1)$=1, where $\pa I_{t}(e^{i\theta})= \pa I(e^{i\theta},t)$;
\item[(iii)] Except for a finite subset of $t\in [0,1]$ not including $0$ or $1$, $\pa I_{t}:\bbS^{1}\hookrightarrow P$
      is transversal to $\pa D_{b}$.
\end{enumerate}

Let $p= \pa I(\bbS^{1},1)\cap \pa D_{b}$ be the unique point of intersection between $\pa I_{1}$ and
$\pa D_{b}$.  Let $b_{1}\in \bbS^{2}$ be an element of $\bbS^{2}$ such that $p\in \pa D_{b_1}$
and such that $\pa D_{b_1}$ is transversal to $\pa D_{b}$.  Assume also that $T_{p}\pa D_{b} \oplus
T_{p}\pa D$ and $T_{p}\pa D_{b}\oplus T_{p}\pa D_{b_{1}}$ have the same orientation.  Choosing $b_{1}$
close enough to $b$, we can also assume that $\pa D_{b}\cap \pa D_{b_1}={p}$.  This is because via 
the twistor correspondence, the set $\pa D_{b}\cap \pa D_{b_1}$ identifies with the set of geodesics
passing through $b$ and $b_{1}$ and if $b_{1}$ is sufficiently close to $b$, there is only one such 
geodesic by lemma~\ref{geod.1}.  

\begin{figure}[h]
\resizebox{5cm}{!}{
\includegraphics[1.5in,0.5in][6in,4.5in]{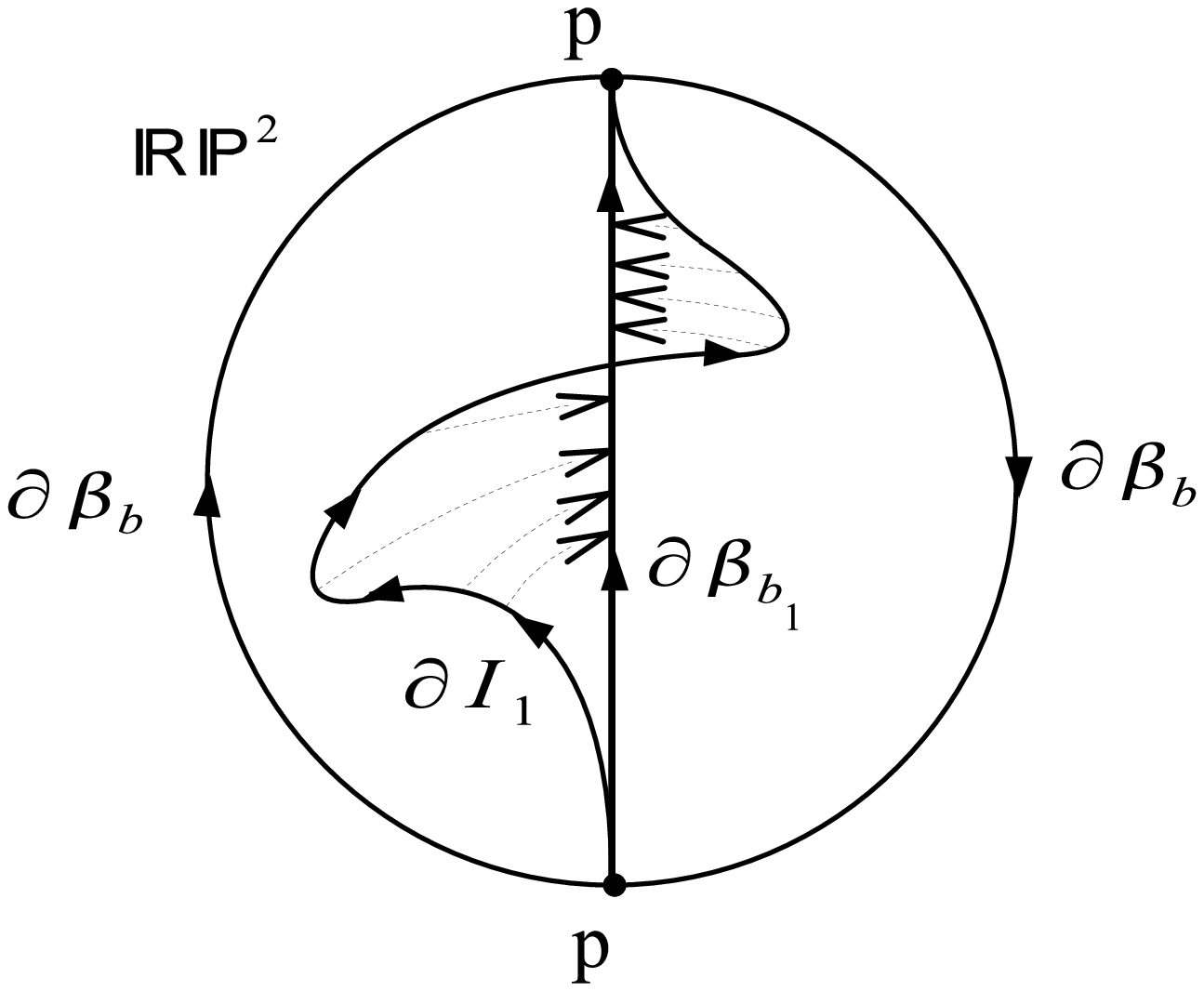}}
	\caption{}
	\label{figure4}
\end{figure}

  Hence, both $\pa I_{1}$ and $\pa D_{b_{1}}$ are embedded circles
in $(P\setminus \pa D_{b})\cup \{p\}$. This is a contractible space, as well as $P\setminus \pa D_{b}$,
which is diffeomorphic to an open disk in $\bbR^{2}$.  Thus, it is possible\footnote{This is intuitively not surprising and can be justified rigorously by invoking the $h$-principle, for instance
theorem 19.4.1, p.176 in \cite{Eliashberg-Mishachev}.} to change the isotopy of 
embeddings $\pa I$ so that beside satisfying the three properties of lemma~\ref{uzs.4}, it is also
such that $\pa I_{1} = \pa \beta_{b_{1}}$ (see figure~\ref{figure4}, where $\rp{2}$  is represented as a disk whose boundary is given by two copies of $\pa D_{b}$ identified via the antipodal map).  By our choice of $b_{1}\in \bbS^{2}$,
we can choose the orientation on $\bbS^{1}$ such that $\pa I_{0}: \bbS^{1}\to \pa D$ and 
$\pa I_{1}: \bbS^{1}\to \pa D_{b_{1}}$ are orientation preserving where $\pa D$ and $\pa D_{b_{1}}$ 
have their orientations induced by the one of $D$ and $D_{b_{1}}$.

If $\theta$ represents the usual angular coordinate on $\bbS^{1}\subset \bbC$, let 
$\frac{\pa}{\pa \theta}\in \CI(\bbS^{1},T\bbS^{1})$ be the corresponding vector field on 
$\bbS^{1}$. Denote by 
\begin{equation}
  \begin{array}{llcl}
       u_{J}: & TP & \to & NP \\
              & v & \mapsto & Jv
  \end{array}
\label{ht.1}\end{equation}
the isomorphism of vector bundles induced by the standard complex 
structure $J$ on $\cp{2}$.  Since $\pa I_{t}$ is an isotopy of embeddings,
\begin{equation}
      u_{J}\left( (\pa I_{t})_{*} \left( \frac{\pa}{\pa \theta}\right) \right)
\label{ht.2}\end{equation}
is a non-vanishing section of $\left. NP\right|_{\pa I_{t} (\bbS^{1})}$, so it defines a section
of the unit vector bundle $\unb P$, this space being canonically identified with $\pa \tcp{2}$.  Thus, the section \eqref{ht.2} defines a lift of $\pa I$ to 
$\pa \tcp{2}$ 
\[
   \pa\tI: \bbS^{1}\times [0,1]\to \pa \tcp{2}, \quad \beta\circ \pa \tI= \pa I.
\]
If $\tD$ and $\tD_{b_{1}}$ are the lifts of $D$ and $D_{b_{1}}$, then by our choices of orientation,
we have that 
\[
   \pa \tI_{0}: \bbS^{1} \to \pa \tD, \quad \pa \tI_{1}: \bbS^{1}\to \pa \tD_{b_{1}}
\] 
are orientation preserving diffeomorphisms.  

Now, working locally in a neighborhood of $P$ in $\cp{2}$, there is no problem in extending the 
isotopy $\pa \tI$ into a smooth homotopy
\[
          \tI: \bbD\times [0,1]\to \tcp{2}, \quad \tI_{t}(e^{i\theta})= \tI(e^{i\theta},t),
          \quad t\in [0,1]
\]
such that $\tI_{0}: \bbD\to \tD$ is a biholomorphism and the restriction of 
$\tI$ to $\pa \bbD\times [0,1]$ is precisely $\pa \tI$.  Let us denote by
\[
             I:= \beta\circ \tI: \bbD\times [0,1]\to \cp{2}
\]
the corresponding map in $\cp{2}$, its restriction to $\pa \bbD\times [0,1]$ being 
precisely $\pa I$.  

\begin{lemma}
We can choose $\tI$ so that there exists $\epsilon>0$ such that 
$\tI(\bbD\times [0,\epsilon])\cap \pa\tD_{b}=\emptyset$ and for all $t>\epsilon$, $I_{t}= \beta\circ \tI_{t}$ is
orthogonal to $P$ with respect to the Fubini-Study metric.
\label{ort.1}\end{lemma}
\begin{proof}
Let $\rho\in \CI([0,\infty))$ be a non-decreasing function such that 
\begin{equation*}
\begin{gathered}
     r<\frac{1}{4} \quad \Longrightarrow \quad \rho(r)=r,  \\
     r> \frac{3}{4} \quad \Longrightarrow \quad \rho(r)=1.
\end{gathered}
\end{equation*}
Consider the smooth homotopy $[0,1]\ni t \mapsto \varphi_{t}$ with $\varphi_{t}(0)=0$ and
\[
         \varphi_{t}(z)= \left(\rho(t)\rho(|z|)+ (1-\rho(t))|z|\right)\frac{z}{|z|}, \quad z\in \bbD^{2}\setminus\{0\}.
\]
Given $\epsilon>0$ such that $\tI(\bbD^{2}\times [0,\epsilon])\cap \tD_{b}=\emptyset$, we can then
consider instead of $\tI$ the smooth homotopy
\[
       \hI= \left\{ \begin{array}{ll}
                     \tI\circ \varphi_{\frac{t}{\epsilon}}, & t\le \epsilon, \\
                     \tI\circ\varphi_{1}, & t\ge \epsilon.
                    \end{array}
       \right.
\]
It is such that $\pa\hI=\pa\tI$ and $\hI_{0}=i$.  Now, let $c:\pa\tcp{2}\times[0,1)_{u}\to \tcp{2}$ be
a collar neighborhood of $\pa\tcp{2}$ induced by the exponential map at $P$ of the Fubini-Study
metric.  Let $\xi\in \CI(\tcp{2},T\tcp{2})$ be a smooth vector field with
$c_{*}\frac{\pa}{\pa u}= \xi$ in a neighborhood of $\pa\tcp{2}$.  Let also
$\psi_{1}, \psi_{2}\in \CI(\bbD^{2}\times [0,1])$ be boundary defining functions for $\bbD^{2}\times \{0\}$
and $(\pa\bbD^{2})\times [0,1]$ respectively, and let $\hI$ evolve according to the ordinary 
differential equation
\[
    \frac{d}{ds} \hI(s;d,t)= \psi_{1}(d,t)\psi_{2}(d,t) \xi(\hI(s;d,t)), \quad
              (d,t)\in \bbD^{2}\times [0,1], \quad \hI(0;d,t):=\hI(d,t).
\]
Then we get the desired homotopy by considering $\hI(s;\cdot,\cdot)$ for $s>0$ small.

\end{proof}

\begin{lemma}
We can choose the smooth homotopy $\tI: \bbD^{2}\times [0,1]\to \tcp{2}$ such that it
is transversal to $\tD_{b}\subset \tcp{2}$ and satisfies the property of lemma~\ref{ort.1}, still demanding that $\tI_{0}: \bbD^{2}\to 
\tD$ is a biholomorphism and that 
\[
     \pa I: \bbS^{1}\times [0,1] \to P
\]
satisfies the three properties of lemma~\ref{uzs.4} and is such that 
$\pa\tI_{1}: \bbS^{1}\to \pa \tD_{b_{1}}$ is a diffeomorphism.
\label{ht.3}\end{lemma} 
\begin{proof}
Let us first investigate the transversality of $\pa\tI$ with $\tD_{b}$ at the boundary.  For $t\in[0,1]$,
a point
$\tp\in\pa\tI_{t}(\bbS^{1})\cap \pa \tD_{b}$ corresponds to a point $p:= \beta(\tp)\in P$ where 
$\pa I_{t}$ is not transversal to $\pa D_{b}$.  Notice however that the converse is not always true.  
At a point $p\in \pa I_{t}(\bbS^{1})\cap \pa D_{b}$ where $\pa I_{t}$ is not transversal to $\pa D_{b}$,
there exists $\tp\in \pa\tI_{t}(\bbS^{1})\cap \pa \tD_{b}$ such that 
$\beta(\tp)=p$ if and only if the orientations of $\pa I_{t}$ and $\pa D_{b}$ at $p$ are the same (so that they give the same lift of $p$ in $\pa\tcp{2}$).  
In any case, we have that
\[
         \pa \tI_{t_{j}}(\bbS^{1})\cap \pa \tD_{b} \ne \emptyset
\]
for a finite number of values of $t$ in $[0,1]$, say $t_{1},\ldots, t_{k}\in [0,1]$.  Let us say
that  $\tp_{1},\ldots,\tp_{k}\in \pa\tI(\bbS^{1}\times [0,1])\cap \pa \tD_{b}$ is
an exhaustive list of points and let $(s_{1},t_{1}),\ldots, (s_{k},t_{k})\in \bbS^{1}\times [0,1]$ be
the unique points such that $\pa \tI_{j}(s_{j},t_{j})=\tp_{j}$.
In fact, by the second property in lemma~\ref{uzs.4}, each non-empty intersection
$\pa \tI_{t_{j}}(\bbS^{1})\cap \pa \tD_{b}$ consists of exactly one point, namely
$\tp_{j}$.  
For each such point, let $p_{j}:= \beta(\tp_{j})$ be the corresponding point in $P$.  To determine
whether or not $\pa \tI$ is transversal to $\pa\tD_{b}$ at $\tp_{j}$, we need to go back to 
remark~\ref{lm.1} about our choice of local model for the isotopy $\pa I$ near the point $p_{j}$.  In
the coordinates $(x,y)$ of the local model of remark~\ref{lm.1}, the point $p_{j}$ corresponds
to the origin.  If $\theta= \arctan\left( \frac{y}{x}\right)$ corresponds to the angular coordinates of
the fibres of $\bbS(TP)$ with the obvious trivialization in the coordinate patch of $(x,y)$, then the 
lifts of $\pa I$ and $\pa D$ to $\bbS(TP)$ at the origin are given respectively by
\begin{equation}
\begin{gathered}
  \pa\hD= \{ (x,y,e^{i\theta}) \in \cU\times \bbS^{1} \quad | \quad y=0, \quad e^{i\theta}=1\}, \\
  \pa \hI_{t}= \{ (x,y,e^{i\theta}) \in \cU\times \bbS^{1} \quad | \quad 
  y= -\epsilon x^{2} +\delta -t, \quad \tan \theta= -2x\epsilon \}.
\end{gathered}
\label{ht.5}\end{equation} 
Thus, at $\hp_{j}=(0,0,1)\in \bbS(TP)$, we have
\[
   T_{\hp_{j}} \pa \hD= \bbR \frac{\pa}{\pa x}, \quad 
   T_{\hp_{j}}\pa \hI = \spa \left\{ \frac{\pa}{\pa y}, \left(\frac{\pa}{\pa x} -2\epsilon \frac{\pa}{\pa \theta} \right) \right\}
\]
so that $\pa \hI$ is transversal to $\pa \hD$.  Identifying $\bbS(TP)$ with $\unb P$ using the 
complex structure $J$ of $\cp{2}$, this means that $\pa \tI$ is transversal to $\pa \tD_{b}$
at $\tp_{j}$.  Since we have this local model for each intersection point $p_{1}, \ldots, p_{k}$,
we conclude that $\pa \tI: \bbS^{1}\times [0,1]\to \pa \tcp{2}$ is transversal to 
$\pa \tD_{b}$ in $\pa \tcp{2}$.

It suffices then to apply the extension theorem~\ref{pre.5} for transversal maps.  Strictly speaking, we would need $\tD_{b}$ to be 
a closed manifold (compact with no boundary), but this can be remedied by considering its double $\widetilde{\bbS}^{2}_{b}$
in the double $\tM$ of $\tcp{2}$.  Defining the double of $\tcp{2}$ in a suitable way (\eg taking
$\tM:= \bbP T_{\bbC}\bbS^{2}$ in \eqref{vvv.1} via the twistor correspondence of LeBrun and Mason), we 
can assume $\widetilde{\bbS}^{2}_{b}\subset \tM$ is $\CI$-embedded in $\tM$.   Then we can first apply the 
extension theorem to make $\tI_{1}$ transversal to $\tD_{b}$ keeping $\pa\tI_{1}= \pa \tD_{b_{1}}$ 
fixed.  But then, clearly the fact that $\pa \tI$ is transversal to $\pa \tD_{b}$ in $\pa \tcp{2}$ 
implies $\pa \tI$ is transversal to $\tD_{b}$ in $\tcp{2}$.  Hence, we can apply the extension
theorem a second time to make $\tI$ transversal, $\tI_{0}$, $\tI_{1}$ and 
$\pa\tI$ being fixed.  In these two applications of the extension theorems, we can
of course preserve the property of lemma~\ref{ort.1}.
\end{proof}

The purpose of the homotopy $\tI$ is to deform $\tD$ so that its boundary coincides with the 
boundary of one of the holomorphic disks of the complete family $\phi: \tcp{2} \to \bbS^{2}$.  This
is because we want to be able to blow-down the boundary of $\tD$ using the blow-down map \eqref{bd.6}.  
However, in order to get a contradiction, we need the oriented intersection number between 
$\tI_{1}$  and $\tD_{b}$ to be positive.  Even though
the homotopy $\tI$ is transversal to $\tD_{b}$, it does not keep track of the intersection number 
since $\tD_{b}$ has a boundary.  Thus, knowing that the oriented intersection number between $\tD$ and $\tD_{b}$ is 
positive does not insure us that the oriented intersection number between $\tI_{1}$ and $\tD_{b}$ is positive.
What saves us is our careful choice of the homotopy $\tI$.

Let $g$ be the Fubini-Study metric on $\cp{2}$.  Then the exponential map induces a tubular
neighborhood
\[
     \nu: \cU\hookrightarrow \cp{2}
\] 
of $P$ in $\cp{2}$ where $\cU\subset \nb P$ is a sufficiently small open set containing the 
zero section, the map $\nu$ identifying the zero section of $\nb P$ with $P\subset \cp{2}$.  
Without loss of generality, we can assume that $\cU$ is invariant under the involution
\[
           \begin{array}{llcl}
            r:& \nb P &\to &\nb P \\
              & (p,n) & \mapsto & (p,-n) 
           \end{array} \quad, n\in \nb_{p}P, \; p\in P.
\] 
Thus, there is a corresponding involution
\[
    \rho: \cU\to \cU
\]
in the tubular neighborhood $\cU \subset \cp{2}$.  Let $\tnb{P}$ be the Melrose's blow-up of the 
zero section of $\nb P$ and let $\alpha_{\nb}: \tnb{P}\to \nb P$ be the blow-down map.  Let
\[
           c: \pa \tnb{P}\times [0,\epsilon)\hookrightarrow \tcU\subset \tnb{P}
\] 
be a collar neighborhood of $\pa \tnb{P}$ in $\tnb{P}$ induced by the tautological 
$[0,\infty)$-bundle over $\pa \tnb{P}\cong \unb P$, whose fibre above
$n\in \unb P$ is the ray spanned by $n\in \nb_{p} P$ where
$p= \pi(n)\in P$ and $\pi: \unb{P}\to P$ is the bundle projection.  Taking $\cU$ smaller if needed, we can assume that for $\tcU:=\alpha^{-1}_{\nb}(\cU)$ there is a 
commutative diagram
\begin{equation*}
\xymatrix{ \tcU \ar[r]\ar[d]^{\alpha} & \tnb{P}\ar[d]^{\alpha_{\nb}}  \\
             \cU \ar[r] & \nb P }  
\end{equation*}
Since the Fubini-Study metric is compatible with the complex structure $J$ of $\cp{2}$, 
\[
        \bD_{b} := (D_{b}\cap \cU) \cup \rho( D_{b}\cap \cU) \subset \cU
\]
has a natural structure of $\mathcal{C}^{1}$-submanifold in $\cU$.  Here, taking $\cU$ small 
enough, we can assume that 
\[
          (D_{b}\cap \cU) \cap \rho( D_{b}\cap \cU)= \pa D_{b} \subset P.
\]
Let $I_{\cU}: I^{-1}(\cU)\to \cU$ be the restriction of the homotopy $I$ to $I^{-1}(\cU)\subset 
\bbD\times [0,1]$.  Let $\bbS^{2}\times [0,1]= (\bbD\cup \bbD_{-})\times [0,1]$ be the 
double of $\bbD\times [0,1]$ induced by some collar neighborhood of $\pa \bbD$ in $\bbD$, 
and let $\cV:= I^{-1}(\cU)\cup (I^{-1}(\cU))_{-}$ be the corresponding double of 
$I^{-1}(\cU)$ on $\bbS^{2}\times [0,1]$.  Then on $\cV$, we can define a homotopy
$I_{\cV}:\cV\to \cU\subset \cp{2}$ by
\begin{equation}
    I_{\cV}(v):= \left\{  \begin{array}{ll}
                          I_{\cU}(v), & v\in I^{-1}(\cU), \\
                          \rho\circ I_{\cU}(v), & v\in (I^{-1}(\cU))_{-} .
    \end{array}
    \right.
\label{ps.2}\end{equation}
By lemma~\ref{ort.1}, the map $I_{\cV}$ is of class $\mathcal{C}^{1}$ on 
$\cV\cap \bbS^{2}\times (\epsilon,1]$.  In fact, directly from the proof of lemma~\ref{ort.1},
one can check that it is even of class $\CI$.

\begin{proposition}
Let $\tI: \bbD\times [0,1]\to \tcp{2}$ be the homotopy of lemma~\ref{ht.3}.  Then
the oriented intersection number between $\tI_{1}: \bbD\to \tcp{2}$ and
$\tD_{b}\subset \tcp{2}$ is positive with intersection points only in $\tcp{2}\setminus \pa
\tcp{2}$.  In fact, it is greater than or equal to the oriented intersection number between
$\tD$ and $\tD_{b}$, which is positive by assumption.  
\label{ps.1}\end{proposition}
\begin{proof}
Let 
\[
    \tS:= \tI^{-1}(\tD_{b})\subset \bbD\times [0,1]
\]
be the preimage of $\tD_{b}$ in $\bbD\times [0,1]$.  Since the map $\tI$ is transversal
to $\tD_{b}$, it is a compact one-dimensional manifold with boundary.  The boundary
of $\tS$ lies in $\bbD\times \{0\}$, $\bbD\times \{1\}$ and $\pa \bbD\times [0,1]$, but not
in the corners $\pa \bbD\times\{0\}$ and $\pa \bbD\times \{1\}$ of the manifold with corners
$\bbD\times [0,1]$.   A choice of orientation of $\bbD\times [0,1]$ together with the orientation
of $\tD_{b}$ impose an orientation on $\tS$.  Clearly, the part of the boundary of $\tS$ contained
in $\bbD\times \{0\}$ gives the intersection points between $\tD$ and $\tD_{b}$, while the part
of the boundary of $\tS$ contained in $\bbD\times \{1\}$ gives the intersection points of 
$\tI_{1}$ with $\tD_{b}$.  If on the other hand $\pa \tS\cap \pa \bbD\times [0,1] =\emptyset$, then
a standard argument of differential topology using the manifold $\tS$ (see figure~\ref{figure5}) shows that the oriented intersection numbers $I(\tI_{0},\tD_{b})$ and 
$I(\tI_{1},\tD_{b})$ are the same, the manifold $\tS$ establishing a cobordism between 
$I(\tI_{0},\tD_{b})$ and $I(\tI_{1},\tD_{b})$.
More generally, if $\pa \tS\cap (\pa \bbD\times [0,1])\ne \emptyset$ (see figure~\ref{figure6}, then the points of 
$\pa S$ lying on $\pa \bbD\times [0,1]$ exactly measure the difference between the oriented 
intersection numbers $I(\tI_{0},\tD_{b})$ and $I(\tI_{1},\tD_{b})$.
\begin{figure}[h]
\begin{minipage}[t]{5.5cm}
\resizebox{4cm}{!}{
\includegraphics[0.5in,1in][5in,5in]{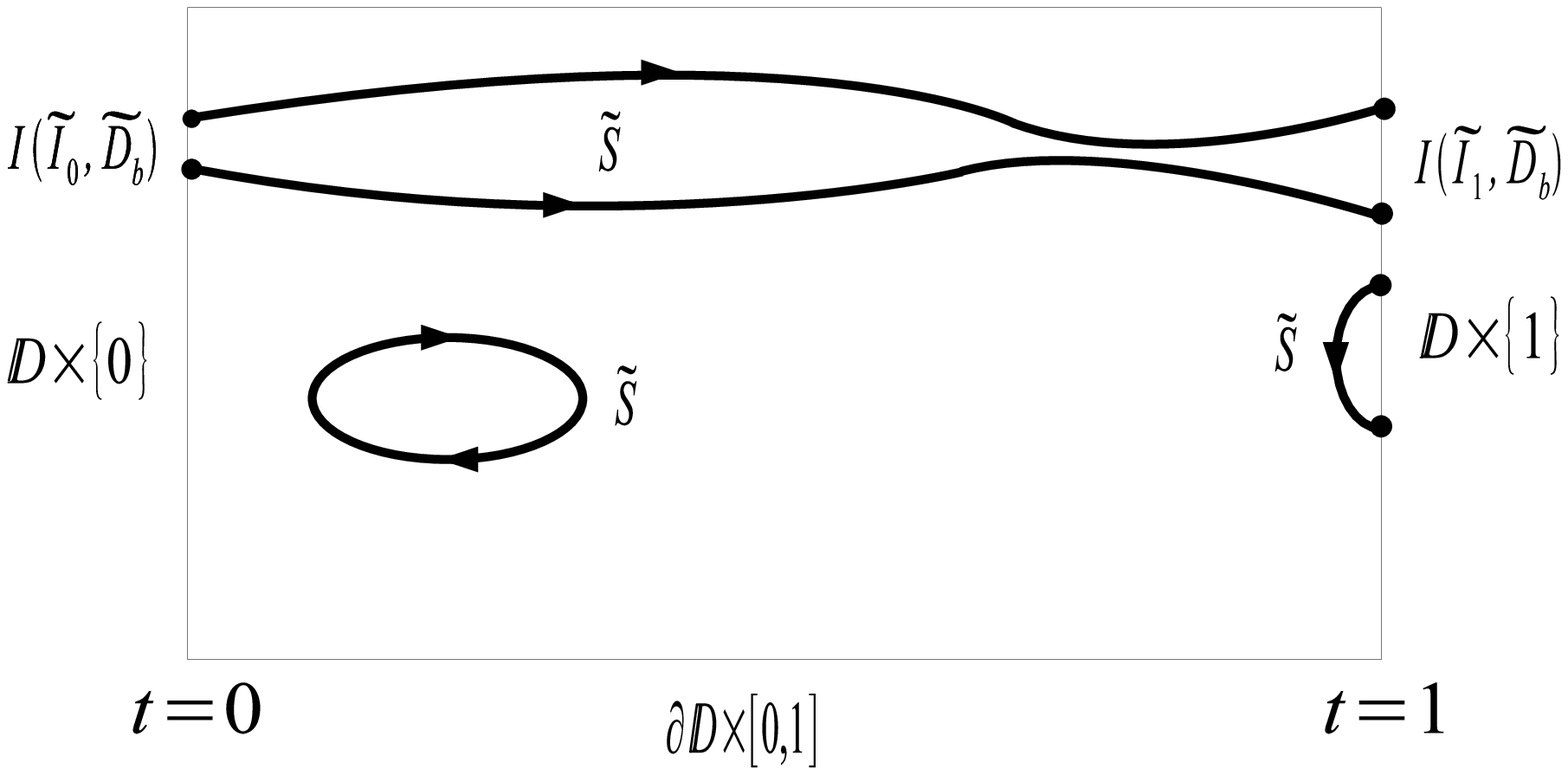}}
	\caption{}
	\label{figure5}
\end{minipage}
\hfill
\begin{minipage}[t]{5.5cm}
\resizebox{4cm}{!}{\includegraphics[0.5in,1in][5in,5in]{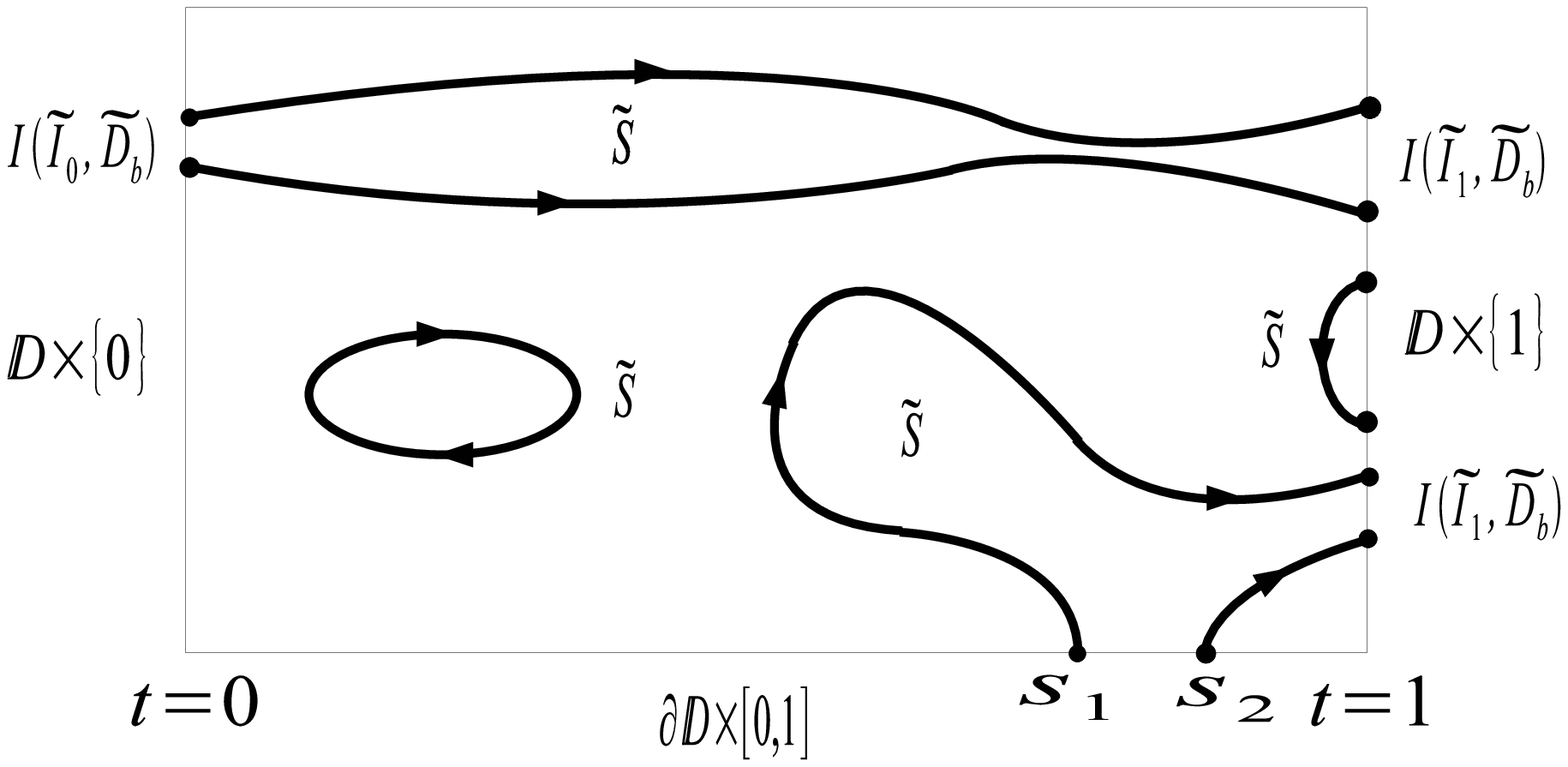}}
	\caption{}
	\label{figure6}
\end{minipage}	
\end{figure}

Let us agree that we choose the orientation of $\bbD\times [0,1]$ so that the induced orientation on $\bbD\times \{0\}$ is such that 
$\tI_{0}:\bbD\times \{0\}\to \tD$ is orientation \textbf{reversing}.  Since the oriented intersection
number between $\tD$ and $\tD_{b}$ is positive, this means that, following the convention of 
\cite{Guillemin-Pollack} p.101 for the orientation of $\tS$, the induced orientation of each point of
$\pa \tS\cap (\bbD\times \{0\})$ is negative.  To prove the proposition, it suffices to show that each point of $\pa \tS\cap (\pa \bbD\times [0,1])$ has a negative induced orientation, since
then
\[
     I(\tI_{1},\tD_{b})= \#[ \pa\tS\cap (\bbD\times \{0\})] + \#[\pa\tS\cap (\pa\bbD\times [0,1])]>0.
\]  
Intuitively, the reason why these points have a negative induced orientation is because, in the
blow-down picture, they come from pairs of intersection points on the boundary at $t=0$ 
with positive intersection numbers.  The details are as follows.

Let $s_{j}=(d_{j},t_{j})\in \pa \tS\cap (\pa \bbD\times [0,1])$, $j\in \{1,\ldots,k\}$ be an
exhaustive list of all the points of the finite set $\pa \tS\cap (\pa\bbD\times [0,1])$.  Near
these points, the map $I_{\cV}$ is of class $\mathcal{C}^{1}$, so it makes sense to discuss 
about transversality there.  To show 
that each of these points has a negative induced orientation from the one of $\tS$, we will
use the blow-down map $\beta: \tcp{2}\to \cp{2}$.  Let $\tp_{j}:=\tI(s_{j})$ be the
corresponding points in $\pa \tcp{2}$ and 
$p_{j}:=\beta\circ \tI(s_{j})= I(s_{j})$ the corresponding points in $P\subset \cp{2}$.  By our
assumptions, $\tI$ is transversal to $\tD_{b}$ at $\tp_{j}$.  On the other hand, this
non-empty intersection indicates that $I=\beta\circ \tI$ is \textbf{not}
transversal to $D_{b}$ at $p_{j}$.  This means that at the point $p_{j}$, the $\mathcal{C}^{1}$-map
$I_{\cV}$ of \eqref{ps.2} is not transversal to $\bD_{b}\subset \cU$.  By the construction
of the homotopy $I$, for each point $s_{j}=(d_{j},t_{j})$, there are two arcs
$(d_{j}'(t),t)$ and $(d_{j}''(t),t)$, where $d_{j}'(t),d_{j}''(t)\in \pa\bbD$ and $0\le t\le t_{j}$, that are mapped by $I$ into $\pa D_{b}$.  These two arcs corresponds
to a pair of intersections points on the boundary that are annihilated via the homotopy
$\pa I$. 

As pointed out in the beginning of the proof of lemma~\ref{ht.3}, $s_{1},\ldots,s_{k}$ is not
an exhaustive list of the values where $I$ fails to be transversal to $D_{b}$.  Let $s_{k+1},
\cdots, s_{k+\ell}\in \pa\bbD\times [0,1]$ be the remaining values for which $I$ fails to 
be transversal to $D_{b}$, and so for which $I_{\cV}$ fails to be transversal to $\bD_{b}$.
These extra points are such that $\tI(s_{j})\notin \tD_{b}$ for $j>k$.      

Now, both $I_{\cV}$ and 
$\bD_{b}$ are invariant under the involution $\rho: \cU\to \cU$, which implies that intersection
points in $\cU\setminus P$ come into pairs of the form $(q, \rho(q))$, $q\in \cU\setminus P$.  
The point $p_{j}$ is on $P$, but it comes from the inductive step of lemma~\ref{uzs.4} to annihilate
two intersection points between $\pa I$ and $\pa D_{b}$ in $P$.

Let $B_{s_{j}}$ be a small ball around $s_{j}$ in $\cV$ for $j\in \{1,\ldots,k+\ell\}$.  By the transversality theorem~\ref{pre.3}, we can
assume $\left. I_{\cV}\right|_{\pa B_{s_{j}}}$ is transversal to $\bD_{b}$.  This is because near $p_{j}$, the map $I_{\cV}$ is transversal to $\bD_{b}$ except at $p$.  Thus, by the extension
theorem~\ref{pre.5}, we can jiggle $I_{\cV}$ inside $B_{s_{j}}$ keeping it
fixed outside $B_{s_{j}}$ so that $\left. I_{\cV}\right|_{B_{s_{j}}}$ becomes transversal to 
$\bD_{b}$.  Notice that this deformation does \textbf{not} have a counterpart for $\tI$.
\begin{figure}[h]
\begin{minipage}[t]{5.5cm}
\resizebox{4cm}{!}{
\includegraphics[0.5in,1in][5in,5in]{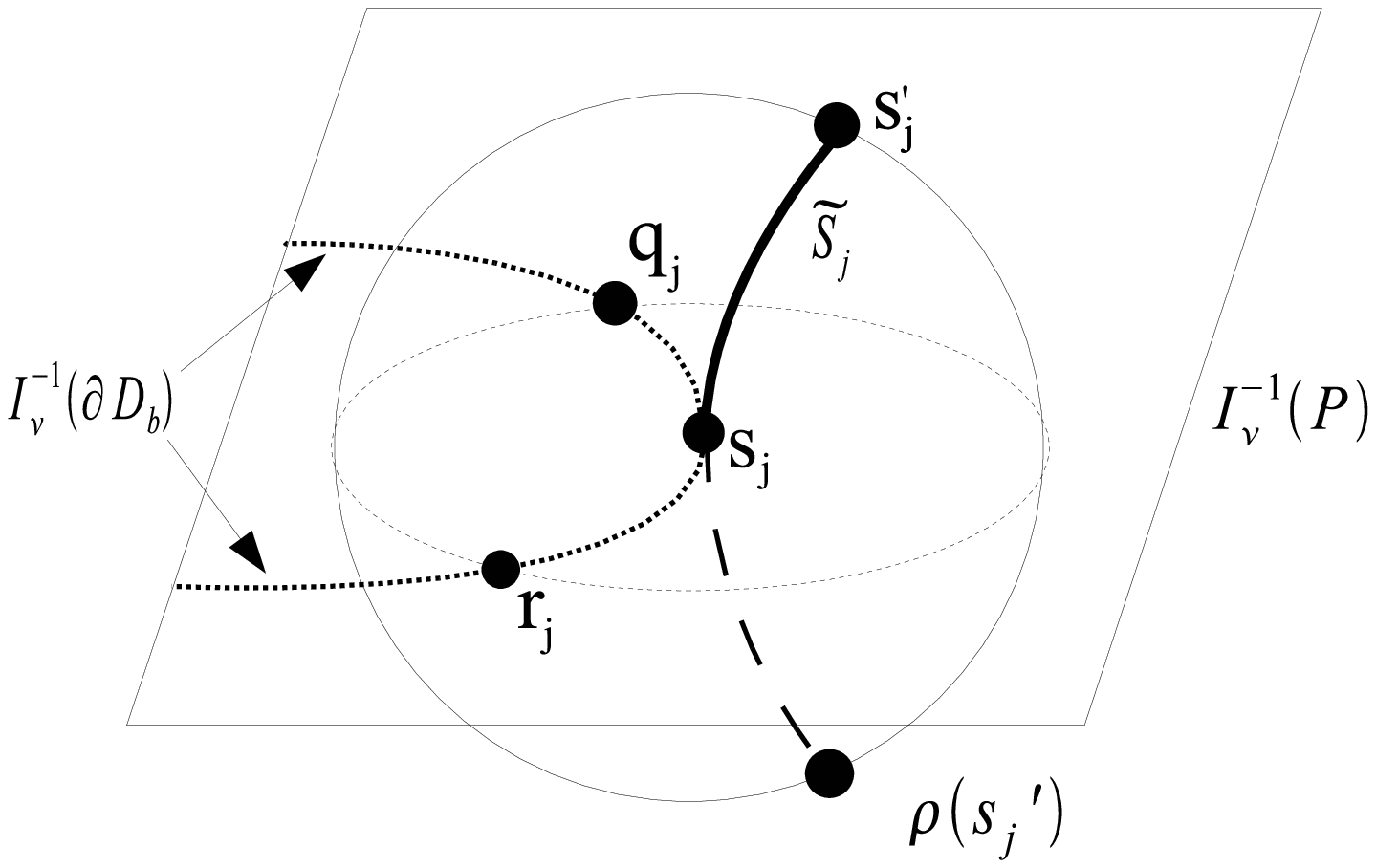}}
	\caption{Before jiggling}
	\label{figure7}
\end{minipage}
\hfill
\begin{minipage}[t]{5.5cm}
\resizebox{4cm}{!}{\includegraphics[0.5in,1in][5in,5in]{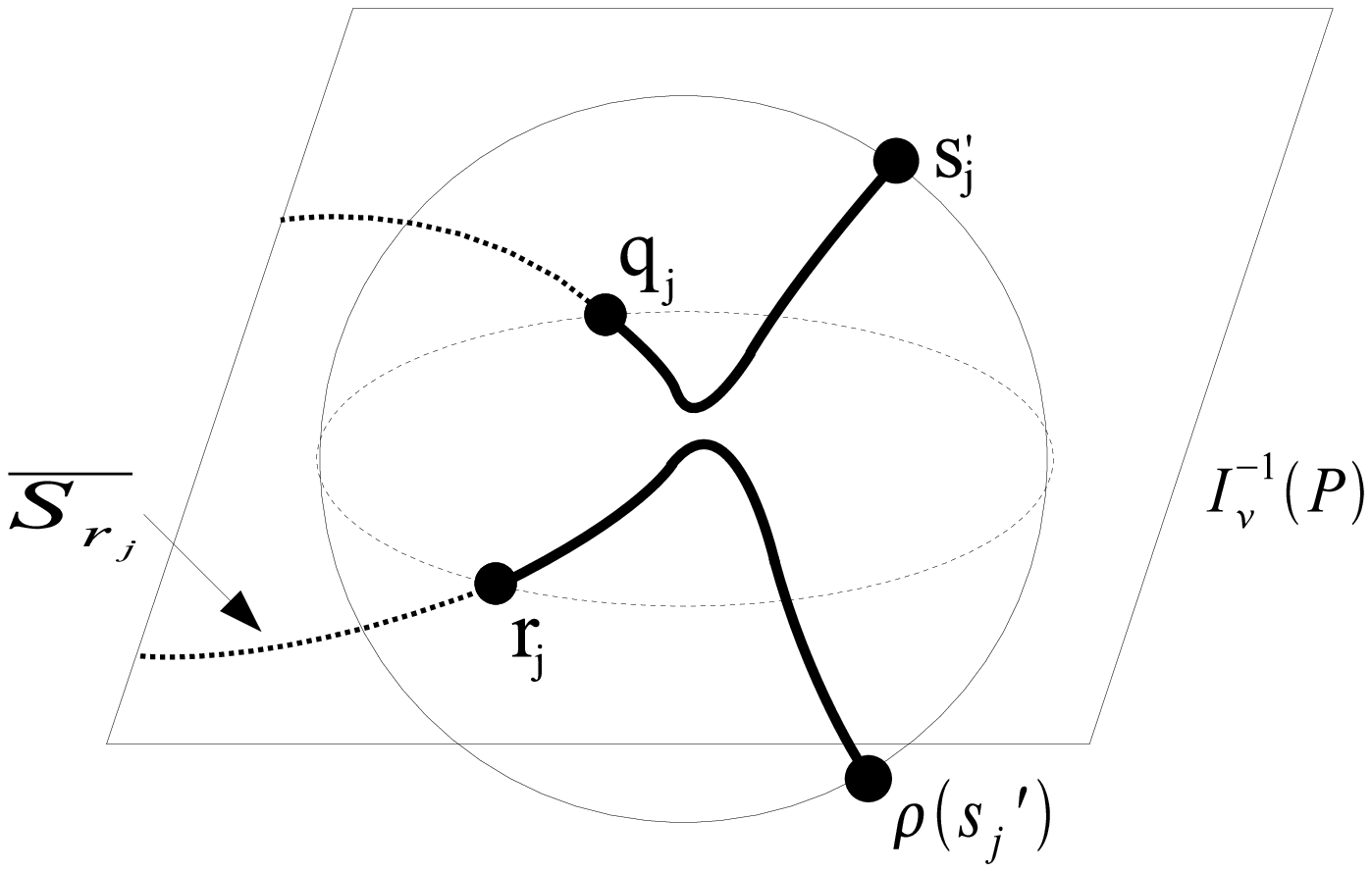}}
	\caption{After jiggling}
	\label{figure8}
\end{minipage}	
\end{figure}

But before doing that, consider the oriented one-manifold with boundary
\[
    \tS_{j}:= \tS\cap B_{s_{j}}, \pa \tS_{j}= \{s_{j}\}\cup \{s_{j}'\}, \quad j\in\{1,\ldots,k\},
\]   
with $s_{j}'\in \pa B_{s_{j}}$ (see figure~\ref{figure7}).  Then the orientation induced by $\tS_{j}$ to $s_{j}$ is the same as
the one induced by $\tS$.  Moreover, the orientation induced by $\tS_{j}$ on $s_{j}'$ is opposite
to the one of $s_{j}$.  On the other hand, when we jiggle $I_{\cV}$ inside each ball
$B_{s_{j}}$ for $j\in\{1,\ldots,k+\ell\}$ to get a new homotopy $\bI_{\cV}$ transversal to $\bD_{b}$, we get a new oriented
one-manifold (see figure~\ref{figure8})
\[
               \bS:= \bI_{\cV}^{-1}(\bD_{b})\subset \cV.
\]
The manifold $\bS\cap (\bbD\times [0,1])$ is the same as $\tS$ away from 
$(\pa\bbD\times [0,1])\cup \cW$ where $\cW= \cup_{j=1}^{k+\ell} B_{s_{j}}$.  Moreover,
$\bS\setminus \cW$ contains additional points that are mapped transversally to 
$\pa D_{b}$ (in $P$) by $\pa I$. 
The manifold $\bS$ is  such that for each $j\in\{1,\ldots, k\}$,
\[
                 \pa (\bS\cap B_{s_{j}})=\{r_{j},q_{j}, s_{j}', s_{j}''\}
\]
consists of exactly four points, where $r_{j},q_{j}\in \pa \bbD\times [0,1]\subset \cV$ and
$\bI_{\cV}(s_{j}'')= \rho( \bI_{\cV}(s_{j}'))$.  Counted with the orientation induced from
$\bS\cap B_{s_{j}}$,  these four points give zero since they are the boundary of the oriented
one-manifold .  Moreover, $s_{j}''$ and $s_{j}'$ have the same induced orientation, since the involution
$\rho: \cU\to \cU$ preserves the orientation of $\cU$ and reverses the orientations of 
$\bD_{b}$ and of $\bI_{\cV}$.  Thus, $r_{j}$ and $q_{j}$ have the same orientation and it is 
 opposite to the one of $s_{j}'$ and $s_{j}''$.  Let $\bS_{r_{j}}$ be the connected component of 
 $\bS\setminus (B_{s_{j}}\cap \bS)$ containing $r_{j}$.  By construction, 
 $\bS_{r_{j}}$ is contained in  $\pa D\times [0,1] \subset \cV$.  Moreover, it is an oriented 
 one-manifold with boundary such that
 \[
          \pa \bS_{r_{j}}= \{r_{j}\}\cup \{\rho_{j}\}, \quad \rho_{j}\in \pa\bbD\times \{0\}.
 \]  
 With our choice of orientation, the orientation of $\rho_{j}$ induced from the one of 
 $\bS_{r_{j}}$ is negative since it corresponds to an intersection point between $D_{b}$ and 
 $I_{0}$, which are holomorphic disks.  Thus, when induced by $\bS_{r_{j}}$, the orientation
 of $r_{j}$ is positive, which means it is negative when induced by $\bI_{\cV}^{-1}(\bD_{b})\cap B_{s_{j}}$.  Thus,
 we conclude that $s_{j}'$ has a positive orientation so that $s_{j}$ has a negative 
 orientation.  Therefore, each point of $\pa S\cap (\pa \bbD\times [0,1])$ has a negative
 induced orientation, which implies that the intersection number between $\tI_{1}$ and 
 $\tD_{b}$ is positive, in fact greater than or equal to $I(\tI_{0},\tD_{b})$.

\end{proof}

Let $\alpha: \tcp{2}\to Y$ be the blow-down map of \eqref{bd.6} associated to the complete family
$\phi: \tcp{2}\to \bbS^{2}$.  By remark~\ref{orientation}, $Y$ comes with a natural orientation.  Since $\pa\tI_{1}:\pa\bbD\to \pa \tD_{b_{1}}$ is a diffeomorphism, there is a 
$\mathcal{C}^{0}$-map
$I_{1}^{\alpha}:\bbS^{2}\to Y$ and a commutative diagram
\begin{equation}
\xymatrix{ \bbD \ar[r]^{\tI_{1}}\ar[d]^{\alpha_{\bbD}} & \tcp{2}\ar[d]^{\alpha}  \\
             \bbS^{2} \ar[r]^{I^{\alpha}_{1}} & Y }  
\label{scc.1}\end{equation} 
where $\alpha_{\bbD}: \bbD\to \bbS^{2}$ is the blow-down map associated to $\bbD$ and the 
disk $\bbD$ is identified with the Melrose's blow-up of $\bbS^{2}$ at a point $s_{0}\in \bbS^{2}$.
Of course, the map $I_{1}^{\alpha}$ is smooth everywhere except possibly at $s_{0}$ where it is
at least continuous.  Similarly,
there is a commutative diagram
\begin{equation}
\xymatrix{ \tD_{b} \ar[r]^{i_{b}}\ar[d] & \tcp{2}\ar[d]^{\alpha}  \\
             \bbS^{2}_{b} \ar[r]^{i_{b}^{\alpha}} & Y }  
\label{scc.2}\end{equation}
where the vertical map on the left is a blow-down map.  Because the boundaries of 
$\tI_{1}(\bbD)$ and $\tD_{b}$ are disjoint, the oriented intersection number of 
$\tI_{1}$ with $\tD_{b}\subset \tcp{2}$ is the same as the one between 
$I^{\alpha}_{1}$ and $\bbS^{2}_{b}\subset Y$.
Now, staying in the same homotopy class, one can approximate $I^{\alpha}_{1}$ by a smooth map.
We can do so by just changing $I^{\alpha}_{1}$ in a small neighborhood of $s_{0}$ so that the intersection points between the new smooth map and $\bbS^{2}_{b}$ are the same as the one
between $I^{\alpha}_{1}$ and $\bbS^{2}_{b}$. 
  In this smooth version of $I^{\alpha}_{1}$,  there is no boundary and the oriented intersection number is a \textbf{homotopy invariant}.  In fact, it only depends on the homology classes of $I^{\alpha}_{1}$ and
  $\bbS^{2}_{b}$ in $\ho_{2}(Y)$.
  
\begin{lemma}
The map $I^{\alpha}_{1}: \bbS^{2}\to Y$ defines the same homology class in $\ho_{2}(Y)$ as
the one associated to $\bbS^{2}_{b}\subset Y$ and $[D]=[\hcf]$.
\label{scc.3}\end{lemma}  
\begin{proof}
Let us first show that $[D]=[h_{\phi}]$.  Since $\tcp{2}$ is a disk bundle over $\bbS^{2}$, $\ho_{1}(\tcp{2})\cong\ho_{1}(\bbS^{2})\cong \{0\}$, so we see from the long exact sequence in homology associated
to the pair $(\tcp{2},\pa\tcp{2})$ that a generator of $\ho_{2}(\tcp{2},\pa\tcp{2})\cong\bbZ$ is sent
to a generator of $\ho_{1}(\pa\tcp{2})\cong\bbZ_{4}$ under the boundary homomorphism
\[
    \pa: \ho_{2}(\tcp{2},\pa\tcp{2})\to \ho_{1}(\pa\tcp{2}).
\]
From the homotopy $\tI$ of lemma~\ref{ht.3}, we see that
\[
     \pa [\tD]= \pa[\tD_{b_{1}}]=\pa[\hcf]. 
\]
This shows that $[D]=[\hcf]$, since otherwise we would get $\pa[\tD]=-\pa[\hcf]$, a contradiction.

Thus, $\tI_{1}$ and $i_{b_{1}}: \tD_{b_{1}}\hookrightarrow \tcp{2}$ define the same
relative homology class in $\ho_{2}(\tcp{2},\pa \tcp{2})$.  

To see that $I^{\alpha}_{1}$ defines the same homology class as $\bbS^{2}_{b}\subset Y$
in $H_{2}(Y)$, notice that there is a natural map
\begin{equation}
     \alpha_{r}: H_{2}(Y)\to H_{2}(\tcp{2},\pa\tcp{2})
\label{hoeq.1}\end{equation}
obtained by composing the maps $H_{2}(Y)\to H_{2}(Y,\alpha(\pa\tcp{2}))$ and 
$\alpha_{*}^{-1}:H_{2}(Y, \alpha(\pa\tcp{2}))\to H_{2}(\tcp{2},\pa\tcp{2})$.  The
homology classes $[I_{1}^{\alpha}], [\bbS^{2}_{b}]\in H_{2}(Y)$ are such that 
\[
      \alpha_{r}(I_{1}^{\alpha}])=\alpha_{r}([\bbS^{2}_{b}])= [\hcf] \quad \mbox{in}
      \; H_{2}(\tcp{2},\pa\tcp{2}).
\]
To see that $[I^{\alpha}_{1}]= [\bbS^{b}]$ in $H_{2}(Y)$, it suffices to check that the 
map \eqref{hoeq.1} is injective.  This follows from the Mayer-Vietoris sequence obtained
by interpreting $Y$ as the union of $\tcp{2}$ with a disk bundle $\mathcal{D}$ over $\bbS^{2}$
with $\pa\mathcal{D}= \pa\cp{2}$,
\begin{equation}
\cdots\longrightarrow H_{2}(\pa\tcp{2})\longrightarrow H_{2}(Y)\longrightarrow 
H_{2}(\tcp{2},\pa\tcp{2})\oplus H_{2}(\mathcal{D},\pa\mathcal{D})\longrightarrow
  \cdots
\label{hoeq.2}\end{equation}
together with the fact that $H_{2}(\pa\tcp{2})=0$ (which follows from the fact the space $\pa\tcp{2}$ is the lens space
$\bbS^{3}/\bbZ_{4}$, \cf the proof of theorem 4.4 in \cite{Lebrun-Mason}).

\end{proof}

\subsection*{Proof of Theorem~\ref{uzs.1}} 

We are now ready to give the proof of theorem~\ref{uzs.1}.  Let $D\subset \cp{2}$ be a holomorphic 
disk as in the statement of the theorem, so such that $\pa D\subset P$,
$D\setminus \pa D \subset P$ and so that its relative homology class $[D]$ is
a generator of $\ho_{2}(\cp{2},P)\cong \bbZ$.  Suppose for a contradiction that 
$D$ is not one of the holomorphic disks of the complete family.  
Then, by lemma~\ref{uzs.3}, there exists $b\in \bbS^{2}$ such that the holomorphic disks $D$ and
$D_{b}:= \beta(\phi^{-1}(b))$ intersects transversely with a positive oriented intersection number in
the interior.  Let $\alpha:\tcp{2}\to Y$ be the blow-down map of \eqref{bd.6}.  
By lemma~\ref{ht.3}, proposition~\ref{ps.1} and lemma~\ref{scc.3}, there exists a 
smooth map $I^{\alpha}_{1}: \bbS^{2}\to Y$ having a positive intersection number with $\bbS^{2}_{b}=
\alpha(\tD_{b})\subset Y$ and corresponding to the same homology class as $\bbS^{2}_{b}$ in 
$\ho_{2}(Y)$.  But the oriented self-intersection number of $\bbS_{b}^{2}$ in $Y$ is given by
\[
   I(\bbS^{2}_{b}, \bbS^{2}_{b})= I(\bbS_{b}^{2}, \bbS_{b_{1}}^{2}) = 0, \quad b_{1}\in \bbS^{2},
\]
since $\bbS^{2}_{b}\cap \bbS^{2}_{b_{1}}=\emptyset$ for $b_{1}\ne b$.  This contradicts the fact that
the oriented intersection number between $I_{1}^{\alpha}$ and $\bbS^{2}_{b}$ is positive.

 \providecommand{\bysame}{\leavevmode\hbox to3em{\hrulefill}\thinspace}
\providecommand{\MR}{\relax\ifhmode\unskip\space\fi MR }
\providecommand{\MRhref}[2]{%
  \href{http://www.ams.org/mathscinet-getitem?mr=#1}{#2}
}
\providecommand{\href}[2]{#2}


\end{document}